\documentclass[journal,twoside,10pt]{IEEEtran}
\usepackage{}
\usepackage{bbm}
\usepackage{mathrsfs} 

\usepackage{relsize}		
\usepackage{url}	                            
\usepackage{wallpaper}		                    
\usepackage{graphicx}
\usepackage{subfigure}
\usepackage{amsmath}
\usepackage{pbox}
\usepackage{booktabs}
\usepackage{bm}
\usepackage{amsmath,amssymb}
\usepackage{algorithm}
\usepackage{algpseudocode}
\usepackage{epsfig}
\usepackage{cite}
\usepackage{textcomp,booktabs}
\usepackage{colortbl}
\definecolor{mygray1}{gray}{.9}
\definecolor{mygray2}{rgb}{.96,.98,.99}
\usepackage{nomencl}
\usepackage{threeparttable}
\usepackage{indentfirst}
\usepackage{multirow}
\newtheorem{remark}{Remark}

\begin{document}

\title{Network-Constrained Transactive Control for Multi- Microgrids-based Distribution Networks with SOPs}

\author{
        Xiaodong~Yang,~\emph{Member,~IEEE,~}
        Zehao Song,~\emph{}
        Jinyu Wen,~\emph{Member,~IEEE,~}
        Chongbo Xu,~\emph{}
        Qiuwei~Wu,~\emph{Senior Member,~IEEE,~}
        Youbing~Zhang,~\emph{Member,~IEEE,~}
        Menglin Zhang,~\emph{}
        and Shijie~Cheng,~\emph{Life Fellow,~IEEE}

\thanks{X.~Yang, J.~Wen and S.~Cheng are with State Key Laboratory of Advanced Electromagnetic Engineering and Technology, School of Electrical and Electronic Engineering, Huazhong University of Science and Technology, Wuhan, 430074, China (e-mail: yang\_xd90@163.com; jinyu.wen@hust.edu.cn; sjcheng@hust.edu.cn). (\emph{Corresponding author: Jinyu~Wen.})
}
\thanks{Y.~Zhang, Z. Song and C. Xu are with the College of Information Engineering, Zhejiang University of Technology, Hangzhou 310023, China (e-mail: youbingzhang@zjut.edu.cn; zehao\_song\_zjut@163.com; chongbo\_xu@163.com).
}
\thanks{Q. Wu is with the Center for Electric Power and Energy, Department of Electrical Engineering, Technical University of Denmark, 2800 Kongens Lyngby, Denmark, and also with the School of Engineering and Applied Sciences, Harvard University, Cambridge, MA 02138 USA (e-mail: qw@elektro.dtu.dk).
}
\thanks{M. Zhang is with the Center for Electric Power and Energy, Department of Electrical Engineering, Technical University of Denmark, 2800 Kongens Lyngby, Denmark (e-mail: menzh@dtu.dk).
}
\vspace{-10pt}
}

\markboth{Journal of \LaTeX\ Class Files,~Vol.~\MakeLowercase{xx}, No.~\MakeLowercase{x}, November~2020}
{Y\MakeLowercase{ang} \MakeLowercase{\emph{et al.}}: Network-Constrained Transactive Control for Multi- Microgrids-based DN\MakeLowercase{s} with SOP\MakeLowercase{s}}

\maketitle

\begin{abstract}
Different from most transactive control studies only focusing on economic aspect, this paper develops a novel network-constrained transactive control (NTC) framework that can address both economic and secure issues for a multi-microgrids-based distribution network considering uncertainties. In particular, we innovatively integrate a transactive energy market with the novel power-electronics device (i.e., soft open point) based AC power flow regulation technique to improve economic benefits for individual microgrids and meanwhile ensure the security of the entire distribution network.
In this framework, a dynamic two-timescale NTC model consisting of slow-timescale pre-scheduling and real-time scheduling stages is formulated to work against multiple system uncertainties.
Moreover, the original bilevel game problems are transformed into single-level mixed-integer second-order cone programming problems through KKT conditions, duality, linearization and relaxation techniques to avoid iterations of transitional methods, so as to improve computational efficiency.
Finally, numerical simulations on a modified 33-bus {\color{black} test system with 3 MGs} verify the effectiveness of the proposed framework.
\end{abstract}


\section{Introduction}

\IEEEPARstart{D}{istribution} networks (DNs) are undergoing a transition from traditional passive systems to active ones due to the increasing promotion of distributed generations, active flexible resources as well as the information and control technologies \cite{feng2020coalitional}. In this context, with the integrating of high-penetration renewable energy sources (RESs), prosumers and various kinds of energy storages, microgrids (MGs) have been emerged as promising self-managed subsystems in DNs for efficient consumption of localized renewable power \cite{shuai2020online,xu2020MultiEnergy}. Recently, multi-MGs (MMGs) are considered as an emerging network designed to further enhance the benefits of MGs \cite{daneshvar2020twostage}. The operation and reliability of the system can be improved by connecting multiple MGs to make a DN with multi-microgrids. Thus, MMGs will be an important network feature in the future DNs, with the merit of operational cost-savings \cite{yang2019interactive}, transmission losses reduction\cite{lezama2019local}, increased utilization efficiency of flexible resources \cite{cui2020efficient}, accommodating more RESs and resilience enhancement \cite{zhou2020flexible}.

Extensive studies have focused on the energy scheduling of MMG systems, which can be classified into two categories: direct control-based and local energy market-based methods.
The control-based methods \cite{farzin2016enhancing,zhao2018energy} are designed from the point of view of the distribution network operator (DNO) or aggregator to determine direct control commands for all controllable parts by collecting all the information from MGs. Although this kind of methods are easy to be implemented, there still have many drawbacks such as privacy concerns, low scalability and also limits the autonomy of individual MGs.

To tackle them, local energy market-based methods emerge for enabling decentralized cooperation among the autonomous MGs, such as, multiagent-based \cite{jadhav2019novel}, game theory-based \cite{anoh2020energy,zhao2020distributed}, auction theory-based \cite{esfahani2019multiagent} methods, and semi-market-based method under supervision \cite{park2016contribution}.

Compared to the direct control-based methods, LEM-based methods tend to provide market platforms that enable the exchange
of energy between the participants only with limited amount of private information exchanges \cite{jin2020local}, which protect self-interests of autonomous MGs and have good scalability.
Nevertheless, above mentioned works only focused on active power optimization for economic operation of MMG systems, while concerning little about the joint investigation of the economic issues of MMG scheduling and the technical issues of DN operation. That is, the distribution network constraints are ignored {\color{black}in \cite{jadhav2019novel,anoh2020energy,zhao2020distributed,esfahani2019multiagent,park2016contribution}}.

Economic and secure operation issues are two major concerns of the operators of DN with MMGs. In practice, the decision-makings of multiple MGs will have an impact on the DN operation, which would affect the energy scheduling of individual MGs conversely if the network constraints including line/transformer capacity and voltage/current limitations are fully considered \cite{wang2020reconfigurable}. Any non-coordination method may lead to the over-utilized of grid assets and may realize suboptimal performances on voltage profiles and overall operational economy.
Meanwhile, the inherent volatility and intermittency of RESs will lead to frequent fluctuation of feeder power, thus resulting in voltage sharp fluctuation and even violation \cite{ji2019robust}.
In addition, the uneven between power generations and consumptions is further aggravated due to the centralized power consumption in certain nodes, which leads to the power imbalance between feeders \cite{ding2018data}, thus disturbing the power flow and increasing system power losses \cite{li2019optimaloperation}. Consequently, these consequences are required to be addressed when optimal managing the DN with MMGs.

Fortunately, transactive control (TC) is emerging as one of the most promising solutions for respecting all the participants' interests \cite{nizami2020multiagent}, including MG owners and DNO.
Transactive control refers to a set of mechanisms for the coordination of various participants through value exchanges, where the price signals are applied to bridge all the components in the system, and the agreement between the control decisions of different components are determined via transactions \cite{liu2020transactive,yang2020transactive}. 
Liu \emph{et al.} \cite{liu2020transactive} proposed a transactive energy based method for the coordinated operation of networked MGs and DNO with distributionally robust optimization.
Yan \emph{et al.} \cite{yan2020distribution} presented a two-level network-constrained method that guarantees the optimal topology of DN and the transactive energy trading among the MMGs.
However, careful review of these excellent studies reveals that the limitation still exists, i.e., the system uncertainties posed by RESs and load demands are not properly addressed. In this regard, any prediction error will lead to inappropriate control commands to physical components, which is undoubtedly detrimental to economic optimality and may even raise security issues to the DN system.

In addition, earlier network-constrained transactive control studies mainly focused on the MG resources management and market clearing algorithms and only took the network security constraints into consideration, where the active measures to further optimize the DN are ignored. In view of this, \cite{wang2020reconfigurable} and \cite{yan2020distribution} proposed to use network reconfiguration to change the topology of the DN for power flow adjustment, which could reduce power losses and mitigate voltage violation. However, limited by the action frequency of line switches, there still has space for performance improvement.
To improve this, soft open points (SOPs), a novel fully-controlled power-electronics device enabling the flexible connection between feeders, is introduced to enhance the flexibility of DN system operations \cite{ji2019robust}.
The SOPs can realize accurately and continuously active/reactive power regulation of the connected feeders and has rapid response speed \cite{cao2016benefits,li2019optimaloperation}, which have been proved more effective than network reconfiguration in power flow adjustment \cite{wang2015SNOP}. Still, the application of SOPs in transactive energy field has not been sufficiently studied previously.

Motivated by aforementioned facts, this paper aims to develop a network-constrained transactive control (NTC) framework for an MMG-based DN system with SOPs under an uncertain environment. This framework organizes a local transactive energy market for both MG owners and DNO to participate in, and respects the interests and preference of individual MGs, and takes active measures to optimize the DN operation to ensure its security.
To our knowledge, we are among the first to explore the benefits of power-electronics device (i.e., SOP) in the network-constrained coordination between the MMGs and DNO, and also address the system uncertainties.

Salient features of this paper are reflected in Table \ref{TabI-1} through the comparison with existing works, and our major contributions are threefold:

\begin{itemize}
  \item Unlike most TC works only focusing on economic aspect, we tend to develop a NTC framework that can address both economic and secure issues for an MMG-based distribution network simultaneously. The active measure of introducing SOPs is innovatively integrated with the AC optimal power flow technique, which could improve economic benefits to the MGs and meanwhile ensure the security of the DN, and has not been studied before.
  \item A dynamic two-timescale NTC model consisting of slow-timescale pre-scheduling and real-time scheduling stages is formulated to work against multiple system uncertainties.
  \item Most decentralized market-based methods clear the prices or make the decisions via iterative algorithms. In contrast, we transform the bilevel game problem into a single-level MISOCP problem through KKT conditions, duality, linearization and relaxation techniques, which results in a significant improvement in computational efficiency.
\end{itemize}

\begin{table}[!t]
\scriptsize
\newcommand{\tabincell}[2]{\begin{tabular}{@{}#1@{}}#2\end{tabular}}
\caption{Comparative Features of Previous studies}
\label{TabI-1}
\centering
\begin{tabular}{c|c|c|c|c}
\hline
Reference & \tabincell{c}{Transactive\\market} & \tabincell{c}{Network\\constraints} & \tabincell{c}{System \\uncertainty} & \tabincell{c}{Active measures \\of the DN}   \\
\hline
\cite{farzin2016enhancing,zhao2018energy} & $-$ & $-$ & Yes & $-$   \\
\hline
\cite{liu2020transactive} & Yes & Yes & $-$ & $-$    \\
\hline
\cite{wang2020reconfigurable,yan2020distribution} & Yes & Yes & $-$ & \tabincell{c}{Yes (network \\reconfiguration)}    \\
\hline
{Proposed NTC} & Yes & Yes & Yes & Yes (using SOPs)    \\
\hline
\end{tabular}
\vspace{-10pt}
\end{table}

The remainder of this paper is organized as follows: Section \ref{sec2} presents the problem description of the DN with MMGs, and outlines the proposed NTC framework. Section \ref{sec3} introduces the optimization models of the DNO and individual MGs, the problem transformation is also given. Detailed formulations of proposed NTC and solution methodology are presented in Section \ref{sec4}. Section \ref{sec5} provides the numerical simulations and Section \ref{sec6} concludes this paper.

\section{Proposed Transactive Control Framework} \label{sec2}

In this study, we consider a common system architecture for a distribution network with MMGs \cite{zhou2020flexible}. The distribution network is divided into several MG areas, and is connected to the up-stream high voltage (HV) system.
In an MG, RESs (mainly the photovoltaic array (PV) and wind turbine (WT)), energy storage system (ESS), general loads and controllable loads are included. 
Besides, the fully controlled power-electronic devices, i.e., SOPs, are installed in the DN for accurately power flow regulation. Both the DNO and MG owners are considered as individual entities in the system and determines its own operation decisions for its controllable resources.
Each MG owner determines the consumption plan of controllable loads, the charing/discharging power of ESS and the exchange power series with the DN to minimize its operational cost. The exchange power represents the imbalance between the demand and generation inside a microgrid that will be compensated by the distribution system with the payment according to clearing price $\lambda_t$. At the same time, the DN operator determines the clearing pricing with the MGs and the active/reactive power outputs of SOPs to improve its operation quality perspective from both economy and security.

\begin{figure}[!t]
  \centering
  \includegraphics[width=2.5in]{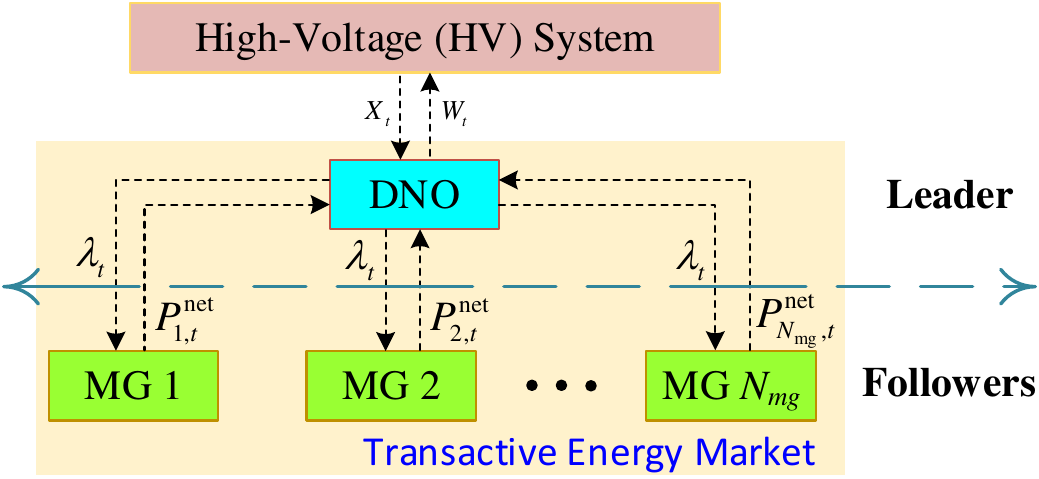}\\
  \caption{Transactive control framework for multi-microgrids-based distribution networks.}
  \label{figfSecII-1}
\vspace{-10pt}
\end{figure}

In view of the above, because the operation of the DNO and each MG is correlated with each other, their operation need to be coordinated in order to achieve the efficient operation of the entire system. In this case, a transactive control based framework is designed for the coordination of the DN operator and MGs, as shown in Fig. \ref{figfSecII-1}.

\begin{figure}[!t]
  \centering
  \includegraphics[width=3.0in]{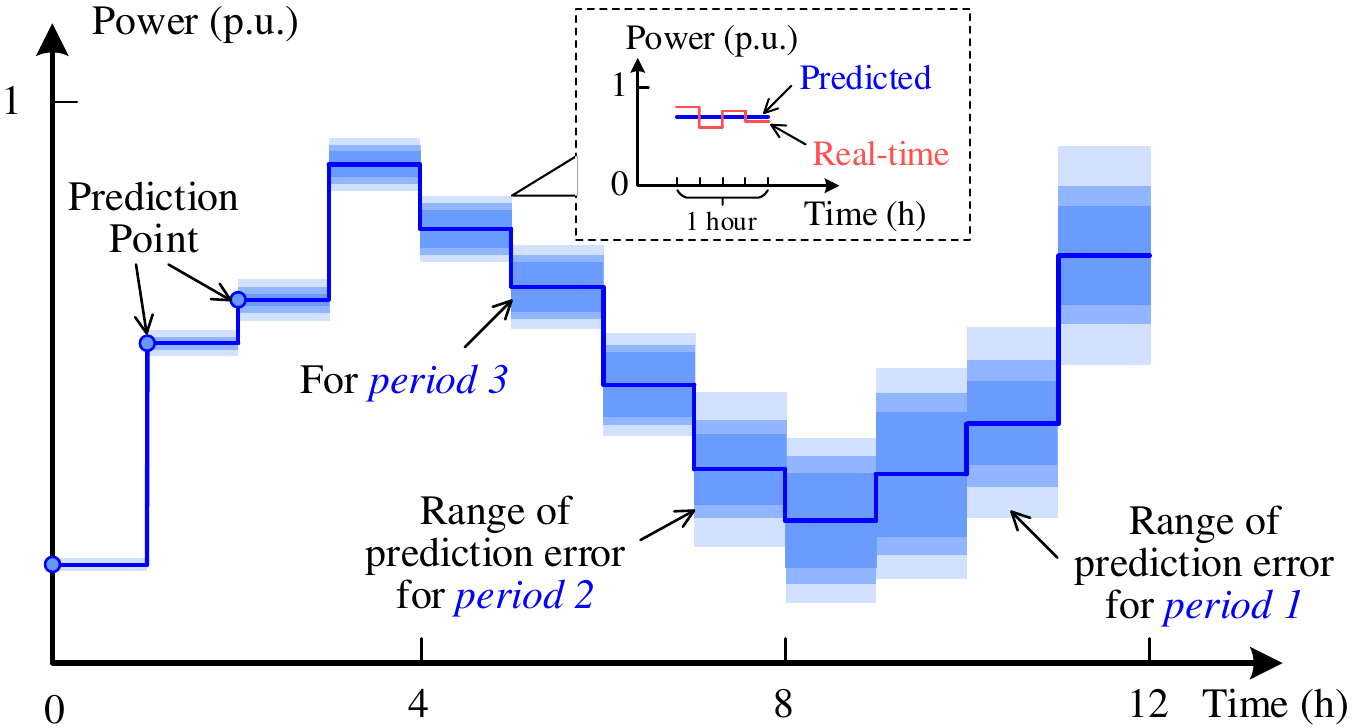}\\
  \caption{Sketch view of prediction errors in different time periods.}
  \label{figfSecII-1.5}
\vspace{-10pt}
\end{figure}

The conclusion of main properties for the focused MMG-based DN is given as follows before designing our approach.
\begin{itemize}
  \item The DNO acts as the leader of the DN system operation and determine the dynamic system prices for the demand and supply in the DN. The MG owners are the followers during the system operation and react to the DNO's decisions. Before the actual real-time implementation, the DNO determines the dynamic system prices and transmits them to all the MGs inside the DN, and the MGs then determine their energy schedules according to the system prices and send their schedules to the DNO. This is a typical iteration process in which the final results are obtained after several iterations.
  \item The high-penetration renewable energy generations are included in the DN, their uncertainties will cause the deviations between the obtained results and expected ones (as shown in Fig. \ref{figfSecII-1.5}), and may also lead to the problems of complex power flow and frequent voltage deviations or even voltage violations.
  \item Lots of previous MMG studies \cite{yang2019interactive,du2020Intelligent} only focus on the active power scheduling while ignoring network topology and active power flow reshaping of the DN. Such non-coordination between MMG and DN operations might lead to over-utilization of DN-side devices, posing the possibility of network congestion and voltage deviation.
\end{itemize}

\begin{figure}[!t]
  \centering
  \includegraphics[width=3.45in]{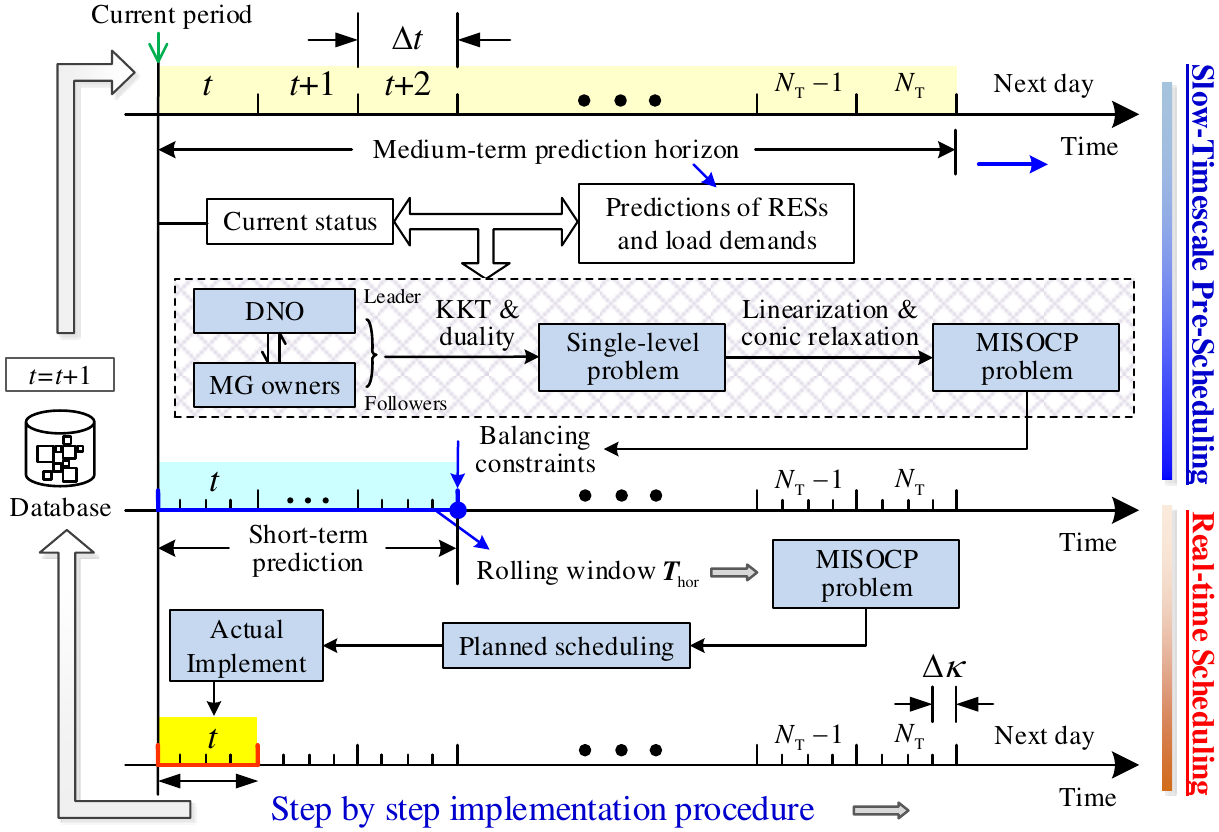}\\
  \caption{Schematic illustration of the proposed method.}
  \label{figfSecII-2}
\vspace{-10pt}
\end{figure}

Considering above, to balance the computational efficiency and the scheduling performances, in the remaining of this paper we focus on providing a mathematical formulation for \emph{two-timescale} network-constrained transactive control method for MMG-based distribution networks, with the coordination of DNO and MGs optimizations. The framework of the proposed NTC method is depicted in Fig. \ref{figfSecII-2}.

{\color{black}
In the proposed NTC framework, the DNO and each MG are considered as individual entities and tend to maximize their own interests. A transacitive energy market is organized by the DNO to coordinate the energy scheduling of the DN and MGs. In particular, instead of fixed prices in the existing models, the power exchange between the DNO and MGs is cleared with dynamic pricing in the transactive market. As such, the energy scheduling of the DN and MGs can be coordinated according to the real-time operation conditions, and the autonomy of individual MGs can be further released and utilized.
}

\section{System Modeling and Model Transformation} \label{sec3}

In this section, the optimization models of DNO and MGs are first formulated, then the formulated bilevel optimization model are transformed into a single level optimization problem through KKT conditions and duality theory.

\subsection{Optimization Model of Distribution Network Operator} \label{sec3.1}

For the DNO, its objectives are twofold: first, to minimize the operational cost; second, to reshape the power flow for power loss reduction and voltage regulation. Also, the operational security of the DN has to be ensured in this process. Therefore, the optimization model of the distribution network operator can be expressed as
\vspace{-2pt}
\begin{align}\label{eqIII-1}
&\min\mathscr{F}_{\mathcal{D}}(\textbf{\emph{x}}) = \mathbbm{a}_{\mathrm{o}}\left(f_{\mathrm{Grid}} + f_{\mathrm{loss}} + f_{\mathrm{sw}} \!-\! \mathscr{I}_{\mathrm{inc}}\right) + \mathbbm{b}_{\mathrm{v}} \mathcal{F}_{\mathrm{vd}}    \\
&\left\{\!\!
\begin{array}{lll} \label{eqIII-2.5}
f_{\mathrm{Grid}} = \sum\nolimits_{t=1}^{N_T}\left(\frac{X_{t}-W_{t}}{2} |g_{t}| + \frac{X_{t}+W_{t}}{2} g_{t}\right)\Delta t;    \\
f_{\mathrm{loss}} = {\mathscr{C}_{\mathrm{loss}}\!\left(\sum\limits_{t=1}^{N_T}\sum\limits_{ij\in \Omega_{l}}^{}r_{ij} I^2_{t,ij}\Delta t + \sum\limits_{t=1}^{N_T}\sum\limits_{i=1}^{N_N}P^{\mathrm{sop,loss}}_{t,i} \Delta t \right)};  \\
f_{\mathrm{sw}} = \sum\limits_{ij\in \Omega_O}\sum\limits_{t=1}^{N_T}(\mathscr{C}_{\mathrm{tap}} |O_{t,ij}-O_{t-1,ij}|);  \\
\mathscr{I}_{\mathrm{inc}} = \sum\limits_{t=1}^{N_T} \sum\limits_{n=1}^{N_{\mathrm{mg}}}{\lambda_t P^{\mathrm{net}}_{t,n}}\Delta t;   \\
\mathcal{F}_{\mathrm{vd}} = \sum\limits_{t=1}^{N_T} \sum\limits_{i=1}^{N_N} \big|U_{t,i}^2-\widetilde{U}_{\mathrm{ref}}^2\big|.
\end{array} \right.
\end{align}
Eq. (\ref{eqIII-1}) is a linear weighted combination of operational cost and voltage deviation minimization problems, where $\mathbbm{a}_{\mathrm{o}}$ and $\mathbbm{b}_{\mathrm{v}}$ are the weight coefficients which can be determined using subjective weighting methods \cite{li2017coordinated}.
$f_{\mathrm{Grid}}$, $f_{\mathrm{loss}}$, $f_{\mathrm{sw}}$ and $\mathscr{I}_{\mathrm{inc}}$ are the grid cost, network losses cost, adjusting cost of on-load tap changer (OLTC), and the income from MGs, respectively.
$X_{t}$ and $W_{t}$ are the trading prices with HV grid; $g_{t}$ is the net load of DN system;
$\mathscr{C}_{\mathrm{loss}}$ and $\mathscr{C}_{\mathrm{tap}}$ are the cost coefficients associated with power losses and OLTC, respectively; $r_{ij}$ and $x_{ij}$ are the resistance/reactance of branch $ij$; $I_{t,ij}$ is the current of branch $ij$ at period $t$; $P_{t,i}^{\mathrm{sop,loss}}$ is the active power losses caused by SOP at period $t$ (see (\ref{eqIII-14}));
$N_T$, $N_N$, and $N_{\mathrm{mg}}$ are the total numbers of time periods, nodes and MGs inside the DN, respectively; $\Delta t$ is the discrete time interval in slow-timescale pre-scheduling stage; $\Omega_l$ and $\Omega_O$ are the sets of all the branches and the branches with OLTCs, respectively;
$\lambda_t$ is the clearing price in transactive market at period $t$; $P_{t,n}^{\mathrm{net}}$ represents the net load of MG $n$ at period $t$;
$\widetilde{U}_{\mathrm{ref}}$ represents the predetermined reference voltage point.

\emph{1) Network constraints of the DN:}

\begin{figure}[!t]
  \centering
  \includegraphics[width=3.45in]{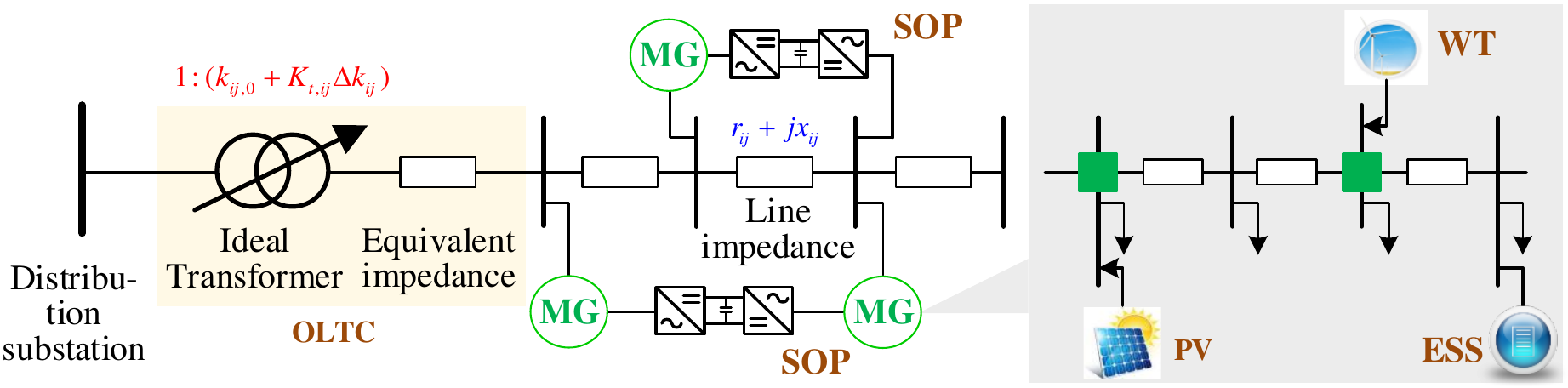}\\
  \caption{Simplified model of MMG-based distribution network with SOPs.}
  \label{figfSecIII-3}
\vspace{-10pt}
\end{figure}

Consider a generic and simplified radial distribution network with SOPs represented in Fig. \ref{figfSecIII-3}, and the widely used \emph{Distflow} branch model is adopted to model the DN \cite{Oikonomou2019deliverable}, as described in the following.
\vspace{-2pt}
\begin{align}\label{eqIII-3}
\sum\nolimits_{ik\in \Omega_l} \!\!\!P_{t,ik} - \sum\nolimits_{ji\in \Omega_l}\!\!(P_{t,ji}-r_{ji}I^2_{t,ji}) = P_{t,i}; ~~~~~~~
\end{align}
\vspace{-12pt}
\begin{align}\label{eqIII-4}
\sum\nolimits_{ik\in \Omega_l} \!\!\!Q_{t,ik} - \sum\nolimits_{ji\in \Omega_l}\!\!(Q_{t,ji}-x_{ji}I^2_{t,ji}) = Q_{t,i}; ~~~~~~
\end{align}
\vspace{-6pt}
\begin{equation}\label{eqIII-5}
U^2_{t,i} - U^2_{t,j} -2(r_{ij}P_{t,ij}+x_{ij}Q_{t,ij}) + (r^2_{ij}+x^2_{ij})I^2_{ij} = 0;
\end{equation}
\vspace{-15pt}
\begin{align}
&I^2_{t,ij}U^2_{t,i} = P^2_{t,ij} + Q^2_{t,ij}; \label{eqIII-6} \\
&P_{t,i} = P_{t,i}^{\mathrm{p}} + P_{t,i}^{\mathrm{w}} + P_{t,i}^{\mathrm{sop}} - P_{t,i}^{\mathrm{L}} + (d_{t,i}^{\mathrm{ess}}-c_{t,i}^{\mathrm{ess}});~~~~~~~  \label{eqIII-7} \\
&Q_{t,i} = Q_{t,i}^{\mathrm{p}} + Q_{t,i}^{\mathrm{w}} + Q_{t,i}^{\mathrm{sop}} - Q_{t,i}^{\mathrm{L}}; \label{eqIII-8}
\end{align}
where $P_{t,i} / Q_{t,i}$ are the sum of active/reactive power injection at node $i$ at period $t$;
$P_{t,ij} / Q_{t,ij}$ are the active/reactive power flow of branch $ij$ at period $t$; $U_{t,i}$ is the voltage magnitude of node $i$ at period $t$; $P_{t,i}^{\mathrm{L}}$ and $Q_{t,i}^{\mathrm{L}}$ are the active/reactive load consumptions of node $i$ at period $t$, respectively; $c_{t,i}^\mathrm{{ess}}$ and $d_{t,i}^\mathrm{{ess}}$ are the charging/discharging power of ESS $i$ at period $t$, respectively;
$P_{t,i}^{\mathrm{p}} / Q_{t,i}^{\mathrm{p}}$ are the active/reactive power injected by PV unit at node $i$, with the relationship of:
\vspace{-3pt}
\begin{equation}\label{eqIII-9}
Q_{t,i}^{\mathrm{p}} = P_{t,i}^{\mathrm{p}} \tan\theta_i^{\mathrm{p}} ~~~~~~~~
\end{equation}
\vspace{-10pt}
\begin{equation}\label{eqIII-10}
\sqrt{(P_{t,i}^{\mathrm{p}})^2+(Q_{t,i}^{\mathrm{p}})^2} \leq S_i^{\mathrm{p}}
\vspace{-2pt}
\end{equation}
where $\theta_i^{\mathrm{p}}$ is the power factor angle of PV $i$, $S_i^{\mathrm{p}}$ is the capacity of PV $i$; and constraints (\ref{eqIII-9})-(\ref{eqIII-10}) are also hold for WTs.

The security constraints of the studied DN are presented as:
\vspace{-3pt}
\begin{equation}\label{eqIII-11}
\underline{U}^2 \leq U^2_{t,i} \leq \overline{U}^2, ~\forall t, i~~~
\end{equation}
\vspace{-10pt}
\begin{equation}\label{eqIII-12}
I^2_{t,ij} \leq \overline{I}^2, ~\forall t, (i,j)\in \Omega_{l}
\vspace{-3pt}
\end{equation}
where $\overline{I}$ is the upper current limit of any branch; $\overline{U}$ and $\underline{U}$ are the upper/lower limits of statutory voltage range, respectively. Constraint (\ref{eqIII-11}) represents the system voltage limits, and the maximum line current capacity is shown in (\ref{eqIII-12}).

\emph{2) SOP operation constraints:}

In this work, the back-to-back voltage source converters (VSCs)-based SOP device is utilized, whose features can be found in \cite{li2017coordinated}.
In this regard, optimization variables of the SOP consist of the active and reactive power outputs of two voltage source converters, thus its flexibility model and constraints are formulated as follows\cite{li2019optimal}:

\emph{i)} Active/reactive power constraints:
\vspace{-3pt}
\begin{equation}\label{eqIII-13}
P_{t,i}^{\mathrm{sop}} + P_{t,j}^{\mathrm{sop}} + P_{t,i}^{\mathrm{sop,loss}} + P_{t,j}^{\mathrm{sop,loss}} = 0, ~\forall t, (i,j)\in \Omega_{\mathrm{S}}
\end{equation}
\vspace{-10pt}
\begin{equation}\label{eqIII-14}
P_{t,i}^{\mathrm{sop,loss}} = A_{i}^{\mathrm{sop}}\sqrt{(P_{t,i}^{\mathrm{sop}})^2 + (Q_{t,i}^{\mathrm{sop}})^2}, ~\forall t
\end{equation}
\vspace{-10pt}
\begin{equation}\label{eqIII-15}
\underline{Q}_{i}^{\mathrm{sop}} \leq Q_{t,i}^{\mathrm{sop}} \leq \overline{Q}_{i}^{\mathrm{sop}}, ~\forall t
\vspace{-3pt}
\end{equation}

\emph{ii)} Capacity constraints:
\vspace{-3pt}
\begin{equation}\label{eqIII-16}
\sqrt{(P_{t,i}^{\mathrm{sop}})^2 + (Q_{t,i}^{\mathrm{sop}})^2}\leq S_{i}^{\mathrm{sop}}, ~\forall t
\vspace{-3pt}
\end{equation}
where $P_{t,i}^{\mathrm{sop}} / Q_{t,i}^{\mathrm{sop}}$ are active/reactive power injection by VSC at node $i$ at period $t$; $\Omega_{\mathrm{S}}$ is the set of branches with SOPs; $S_{i}^{\mathrm{sop}}$ is the capacity of VSC at node $i$; $\underline{Q}_{i}^{\mathrm{sop}}$ and $\overline{Q}_{i}^{\mathrm{sop}}$ are reactive power boundaries of VSC at node $i$.
It should be noted that (\ref{eqIII-13})--(\ref{eqIII-16}) are also hold for node $j|(i,j)\in \Omega_{\mathrm{S}}$.

\emph{3) Constraints of the OLTC:}

The equivalent model of an OLTC is depicted in the left side of Fig. \ref{figfSecIII-3}. The optimization variables of the OLTC is considered as the action series of its tap steps (which is represented as $O_{t,ij}$). Considering this, the operation constraints of the OLTC are formulated as follows.
\vspace{-5pt}
\begin{align}\label{eqIII-17.1}
& U_{t,j} = (o_{ij,0}+O_{t,ij}\Delta o_{ij})U_{t,i}; \\
& \sum\nolimits_{t=1}^{N_T}\big| O_{t,ij}-O_{t-1,ij} \big| \leq \overline{\Delta}^{\mathrm{OLTC}}; \\
& -\overline{O}_{ij} \leq O_{t,ij} \leq \overline{O}_{ij}.  \label{eqIII-17.3}
\end{align}
where $o_{ij,0}$ and $\Delta o_{ij}$ are the initial turn ratio and increment of the OLTC; $\overline{O}_{ij}$ is the total tap steps of the OLTC; $\overline{\Delta}^{\mathrm{OLTC}}$ is the number of allowed actions for the OLTC.

\emph{4) Transactive energy price constraints:}

Following boundary limitation should be kept when clearing the transactive energy prices.
\begin{equation}\label{eqIII-17.5}
\lambda_t^{\mathrm{min}} \leq \lambda_t \leq \lambda_t^{\mathrm{max}}
\end{equation}
where $\lambda_t^{\mathrm{min}}$ and $\lambda_t^{\mathrm{max}}$ are the boundaries of clearing price at period $t$.

\subsection{Optimization Model of Microgrids} \label{sec3.2}

The optimization model of microgrids is responsible for managing the consumption plans of load demands and the usage of ESS in the transactive market organized by DN, to minimize its operational cost. With this objective settings, the optimization model of individual MG $n$ is given as below.
\vspace{-3pt}
\begin{align}\label{eqIII-18}
&\min \mathcal{J}_n \!=\! \sum\limits_{\forall t}{\lambda_t P^{\mathrm{net}}_{t,n}}\Delta t + \sum\limits_{\forall t}\sum\limits_{i\in \mathcal{S}_{n}} \vartheta_n^{\mathrm{Disc}}(P_{t,i}^{\mathrm{L,do}} + P_{t,i}^{\mathrm{L,up}}) \Delta t   \nonumber \\
&+ \sum_{\forall t}\sum_{i\in \mathcal{S}_{n}} (c_{t,i}^\mathrm{{ess}}\eta^\mathrm{{ess,c}}+\frac{d_{t,i}^\mathrm{{ess}}}{\eta^\mathrm{{ess,d}}})\Delta t \mathscr{C}_{\mathrm{deg}}.
\end{align}
where $\vartheta_n^{\mathrm{Disc}}$ is the inconvenience sensitivity coefficient of MG $n$; $P^{\mathrm{L, up}}_{t,i} / P_{t,i}^{\mathrm{L,do}}$ are the increased/decreased load demands for DR, respectively; $\eta^{\mathrm{ess,c}}$ and $\eta^{\mathrm{ess,d}}$ are the power exchange efficiencies; $\mathcal{S}_{n}$ represents the node set of MG $n$; $\mathscr{C}_{\mathrm{deg}}$ is the coefficients concerning ESS degradation.

\emph{1) Operation constraints of ESS:}
The charging/discharing power and the state of charge of ESSs must meet the following constraints to ensure their normal operation.
\begin{align}\label{eqIII-23}
\left\{\!\!
\begin{array}{lll}
0\leq c_{t,i}^{\mathrm{ess}} \leq u_{t,i}^{\mathrm{ess}} c_{i}^{\mathrm{rat}} \\
0\leq d_{t,i}^{\mathrm{ess}} \leq (1-u_{t,i}^{\mathrm{ess}})d_{i}^{\mathrm{rat}}
\end{array} \right. \forall t, i
\end{align}
\begin{equation}\label{eqIII-23.5}
c_{t,i}^{\mathrm{ess}} d_{t,i}^{\mathrm{ess}} = 0
\end{equation}
\begin{equation}\label{eqIII-24}
S_{t,i} = S_{t-1,i}+\frac{c_{t,i}^{\mathrm{ess}}\eta^{\mathrm{ess,c}} - d_{t,i}^{\mathrm{ess}}/\eta^{\mathrm{ess,d}}}{Cap_{i}^{\mathrm{ess}}}\Delta t, ~\forall t, i
\end{equation}
\vspace{-10pt}
\begin{equation}\label{eqIII-25}
S_i^{\min} \leq S_{t,i} \leq S_i^{\max}, ~\forall t,i
\end{equation}
\vspace{-10pt}
\begin{equation}\label{eqIII-25.5}
S_{1,i} = S_{N_T,i}, ~\forall i
\end{equation}
where $S_{t,i}$ is the state of charge (SoC) of ESS $i$ at period $t$; $Cap_{i}^{\mathrm{ess}}$ is the capacity of ESS $i$; $c_{i}^{\mathrm{rat}}$ and $d_{i}^{\mathrm{rat}}$ are the rated charging/discharging power; $S_i^{\min}$ and $S_i^{\max}$ are the SoC boundaries of ESS $i$.

\emph{2) Operation constraints of the demand response (DR) resources:}
To account for DR programs, the following operation constraints are introduced \cite{yang2019interactive}.
\begin{equation}\label{eqIII-26}
\left\{\!\!\!
\begin{array}{ll}
L_{t,i}^{\mathrm{min}}\leq P_{t,i}^{\mathrm{L,up/do}} \leq L_{t,i}^{\mathrm{max}},  &  t\in [t_i^{\mathrm{l}},t_i^{\mathrm{r}}], \forall i \\
P_{t,i}^{\mathrm{L,up/do}} = 0,  &  ~\!\!t \notin [t_i^{\mathrm{l}},t_i^{\mathrm{r}}], \forall i
\end{array} \right.  ~~~~~~~~~~~~
\end{equation}
\vspace{-5pt}
\begin{equation}\label{eqIII-27}
\sum\nolimits_{t} \!P^{\mathrm{L, up}}_{t,i} = \sum\nolimits_{t} \!P_{t,i}^{\mathrm{L,do}}, ~\forall i  ~~~~~~~~~~~~~~~~~~~~~~~~~~~~~~~~
\end{equation}
\vspace{-10pt}
\begin{equation}\label{eqIII-27.5}
P^{\mathrm{L, up}}_{t,i} P_{t,i}^{\mathrm{L,do}} = 0, ~\forall t, i   ~~~~~~~~~~~~~~~~~~~~~~~~~~~~~~~~~~~~~
\end{equation}
\vspace{-10pt}
\begin{equation}\label{eqIII-28}
Q_{t,i}^{\mathrm{L}} = (P^{\mathrm{L,unc}}_{t,i}+P^{\mathrm{L, up}}_{t,i}+P_{t,i}^{\mathrm{L,do}}) \tan(\cos^{-1}(pf_{\mathrm{l}})), ~\forall t, i
\vspace{-3pt}
\end{equation}
where $pf_{\mathrm{l}}$ is the average power factor of load demands; $P^{\mathrm{L,unc}}_{t,i}$ is the uncontrollable load demands; $Q_{t,i}^{\mathrm{L}}$ is the reactive load demands at node $i$ at period $t$, respectively; $[L_{t,i}^{\mathrm{min}},L_{t,i}^{\mathrm{max}}]$ is the range of shiftable loads; $[t_i^{\mathrm{l}},t_i^{\mathrm{r}}]$ is the expected operation time range of shiftable loads at node $i$.

\subsection{Transformation of the Bilevel Model}

As the models formulated in Sections \ref{sec3.1} and \ref{sec3.2}, the combined optimization of DNO and MGs is a typical \emph{bilevel problem} which can be expressed as
\vspace{-2pt}
\begin{align}\label{eqIII-30}
Upper Level: ~&\min ~ \mathscr{F}_{\mathcal{D}}(\textbf{\emph{x}}); \nonumber \\
\mathrm{Subject ~to:} ~&(\ref{eqIII-3})-(\ref{eqIII-17.5});  \nonumber \\
\mathrm{Variables:} ~&\lambda_t, O_{t,ij}, P_{t,i}^{\mathrm{sop}}, Q_{t,i}^{\mathrm{sop}};~~~~~~~~~   \nonumber
\end{align}
\begin{equation}\label{eqIII-31}
\lambda_t \!\downdownarrows ~ \upuparrows \!P^{\mathrm{net}}_{t,n}, \forall n ~~~~~~~~~   \nonumber
\end{equation}
\vspace{-20pt}
\begin{align}\label{eqIII-32}
Lower Level: ~&\min ~\mathcal{J}_n, ~\forall n; \nonumber \\
\mathrm{Subject ~to:} ~&(\ref{eqIII-23})-(\ref{eqIII-28});  \nonumber \\
\mathrm{Variables:} ~&P^{\mathrm{net}}_{t,n}, P^{\mathrm{L, up}}_{t,i}, P_{t,i}^{\mathrm{L,do}}, c_{t,i}^\mathrm{{ess}}, d_{t,i}^\mathrm{{ess}}.   
\end{align}

Generally, the bilevel optimization problem (\ref{eqIII-32}) can be solved via Stackelberg game theory \cite{liu2018energysharing}.
The clearing price $\lambda_t$ is determined by the DNO and serves as an input of the MGs' optimization problem. Based on the received $\lambda_t$, the operation of the MGs are determined by the MGs' optimization, and the power exchange between the DNO and MGs are also determined accordingly. At the same time, the power exchange must satisfy the network constraints and match the solution of the MGs' optimization.
Such a method needs to obtain an equilibrium solution after an iterative process.

In order to efficiently solve the bilevel optimization problem (\ref{eqIII-32}), the KKT conditions of the MGs' optimization are used to transform the original bilevel optimization into a single level optimization problem.

Before transforming the bilevel optimization problem, the complementary relaxation constraints of ESS and DR program need to be discussed. As \cite{li2016sufficient} has proven that the nonlinear constraint (\ref{eqIII-23.5}) can be relaxed and removed from the model, and is exact for the global optimal solution, the proof is omitted here. Similarly, the complementary relaxation constraint {\color{black}(\ref{eqIII-27.5})} concerning the DR can also be removed safely.

On these bases, the leader-followers game problem (\ref{eqIII-32}) can be transformed via KKT conditions of MGs' optimization model and duality theory. (\ref{eqIII-32}) is thus equivalent to the follows.
\begin{align}\label{eqIII-33}
&\min \mathscr{F}_{\mathrm{1equ}} = \left\{\mathbbm{a}_{\mathrm{o}}\left(f_{\mathrm{Grid}} \!+\! f_{\mathrm{loss}} \!+\! f_{\mathrm{sw}} \!-\! \mathscr{I}_{\mathrm{inc}}\right) + \mathbbm{b}_{\mathrm{v}} \mathcal{F}_{\mathrm{vd}} \right.   \nonumber \\
&+\sum_{t} \left\{\sum_{i}[\mu_{1,t,i}^{L}L_{t,i}^{\mathrm{min}}-\mu_{2,t,i}^{L}L_{t,i}^{\mathrm{max}}]   \right. \nonumber \\
&\left. + \sum_{e}(-\mu_{2,t,e}^{\mathrm{ess,c}}c_{e}^{\mathrm{rat}} -\mu_{2,t,e}^{\mathrm{ess,d}}d_{e}^{\mathrm{rat}} +\mu_{1,t,e}^{\mathrm{ess}}S_{e}^{\mathrm{min}} \!-\! \mu_{2,t,e}^{\mathrm{ess}}S_{e}^{\mathrm{max}})\right\}   \nonumber \\
& {\color{black}+ \lambda_{2,e}\sum_{e}S_{1,e}}     \\
&\mathrm{Subject ~to:} ~(\ref{eqIII-3})-(\ref{eqIII-17.5});  \nonumber \\  \label{eqIII-33.1}
& \textbf{0} \leq \bm{\mu} \bot \bm{h(x)}\geq \textbf{0}~(\mathrm{for}~(\ref{eqIII-23}),(\ref{eqIII-24})-(\ref{eqIII-27}),(\ref{eqIII-28}))     \\
&\mathrm{Variables:} ~\lambda_t, O_{t,ij}, P_{t,i}^{\mathrm{sop}}, Q_{t,i}^{\mathrm{sop}}, P^{\mathrm{net}}_{t,n}, P^{\mathrm{L, up}}_{t,i}, P_{t,i}^{\mathrm{L,do}}, c_{t,i}^\mathrm{{ess}}, d_{t,i}^\mathrm{{ess}}.   \nonumber
\end{align}
where $\mu$ and $\lambda$ are the dual variables of inequality and equality constraints in optimization model of MGs, respectively.
It is worth noting that the standard form (\ref{eqIII-33.1}) is available for the transitions of constraints (\ref{eqIII-23}), (\ref{eqIII-24})-(\ref{eqIII-27}) and (\ref{eqIII-28}).

\section{Formulation of Proposed NTC Method} \label{sec4}

On the basis of Sec. III, this section presents the mathematical models for slow-timescale pre-scheduling and real-time scheduling stages of the proposed NTC method to against system uncertainties, and also provides the solution methodology and implementation algorithm.

\subsection{Uncertainty Modeling}

There are a number of uncertainties from load profiles and variable on-site renewables (i.e., PV and WT generators) that could potentially affect the scheduling decisions of the DN and MGs. In this work, we focus on the network-constrained transactive control between the DN and MGs, thus, prediction techniques for uncertain factors are out of scope for this paper.
Therefore, a simple exponential smoothing model \cite{yang2019real} is used to predict short-term renewable outputs and load demands $g_{t,i,r}$ based on the historical data.
The active outputs of the WT generator and PV array, as well as the active load demands can be uniformly expressed as
\vspace{-5pt}
\begin{equation}\label{eqA35}
g_{t,i,r} = \xi\,\overline{g}_{t,i,r} + (1-\xi) \widehat{g}_{t,i,r}, ~r=1,2,3
\end{equation}
\vspace{-10pt}
\begin{equation}\label{eqA36}
\widehat{g}_{t,i,r} = \overline{g}_{t,i,r}(1 + u_{i,r}\cdot \Upsilon_{i,r})~~~~~~~~~~
\end{equation}
\vspace{-15pt}

\noindent where $r=1,2,3$ corresponds to WT, PV, and load demands, respectively, that is $g_{t,i,r}=[P_{t,i}^{\mathrm{w}},P_{t,i}^{\mathrm{p}},P^{\mathrm{L}}_{t,i}]^{\mathrm{T}}$, where $P^{\mathrm{L}}_{t,i}$ is the load demands at node $i$ at period $t$; $\overline{g}_{t,i,r}$ and $\widehat{g}_{t,i,r}$ are the predicted value and its corresponding  stochastic variable; $\xi$ is a predetermined coefficient, such that $0 < \xi < 1$;
$\Upsilon_{i,r}$ is a random number that follows specific normal distribution \cite{jiang2019stochastic}; $u_{i,r}$ is the uncertainty percentage of the RESs outputs or the load demands. 

\begin{remark}
It remains true that the global optimal decision-making of any scheduling is made on the basis of prediction information. In this regard, the outputs of optimization problem (\ref{eqIII-33}) will deviated from expected ones due to the always existed prediction errors caused by the inherent uncertainties of RESs and load demands.
As shown in Fig. \ref{figfSecII-2}, the prediction error increases with the time distance from the prediction point increases. Meanwhile, the RES power and load demands change frequently even within an hour because of their strong volatility and intermittency.
These facts create the exact needs of a multi-timescale TC architecture for an MMG-based DN to against the system uncertainties in a dynamic manner.
\end{remark}

\subsection{Two-Timescale NTC Model}

The proposed NTC method will be modeled from following two sub-models with different timescales.

\emph{1) Slow-Timescale Pre-Scheduling Model}

The pre-scheduling optimization stage of the NTC method is conducted in slow-timescale in a dynamic manner to provide pre-scheduling strategies from the perspective of global optimization in long time horizon, which can provide reference strategy and long-view guidance for the real-time scheduling stage. In this regard, at every period $t$, the slow-timescale pre-scheduling model can be formulated in the following.
\vspace{-2pt}
\begin{align}\label{eqIV-1}
~&\min ~ \mathscr{F}_{\mathrm{1equ}}\big|_{\Delta t}^{t\rightarrow N_T};  \\
\mathrm{Subject ~to:} ~&(\ref{eqIII-3})-(\ref{eqIII-17.5}), ~{\color{black}(\ref{eqIII-33.1})};  \nonumber
\end{align}
where symbol $\diamond\big|_{\Delta t}^{t\rightarrow N_T}$ means that the terms in (\ref{eqIV-1}) (the details can be seen in (\ref{eqIII-2.5}) and (\ref{eqIII-33})) are calculated with period $\Delta t$ from period $t$ to $N_T$.

\emph{2) Real-Time Scheduling Model}

In order to further work against the time-variant and uncertain features of the renewable energy generations and load demands, the real-time scheduling stage is designed nesting the previous stage. With the purpose of improving the control precision and implementary efficiency of the NTC method, this stage is conduced in fast-timescale via a rolling receding manner. The time interval $\Delta t$ is divided into $\zeta$ shorter sampling periods with an interval of $\Delta \kappa$, so the $t$-th $\Delta t$ $\equiv$ $\{(t-1)\zeta+1, \ldots, t\zeta\}\Delta \kappa$. Assuming the real-time scheduling stage is implemented in a short time horizon with $T\Delta t$, that is $\textbf{\emph{T}}_{\mathrm{\textbf{hor}}}=\{(t-1)\zeta+1, \ldots, (t+T-1)\zeta\}$. At each time period $t$, based on the short-term predictions of the upcoming load demand and generations in $\textbf{\emph{T}}_{\mathrm{\textbf{hor}}}$, the Real-Time Scheduling optimization is formulated to determine the real-time scheduling decisions for horizon $\textbf{\emph{T}}_{\mathrm{\textbf{hor}}}$, but only the decisions at current period $t$ is actually implemented on the distribution network.

Considering the above, the real-time scheduling model is thus formulated as follows.
\vspace{-3pt}
\begin{align}
~&\min ~ \mathscr{F}_{\mathrm{1equ}}\big|^{(t-1)\zeta+1\rightarrow(t+T-1)\zeta}_{\Delta \kappa};  \label{eqIV-3} \\
\mathrm{Subject ~to:} ~&(\ref{eqIII-3})-(\ref{eqIII-17.5}),~{\color{black}(\ref{eqIII-33.1})};  \nonumber \\ \label{eqIV-4}
&S_{(t+T-1)\zeta,e}\big|_{\Delta \kappa} = S_{(t+T-1),e}\big|_{\Delta t}, ~\forall e   \\
&\sum\limits_{h\in \textbf{\emph{T}}_{\mathrm{\textbf{hor}}}} \!\!\!(P^{\mathrm{L, up}}_{h,i} - P_{h,i}^{\mathrm{L,do}})\Delta \kappa =  \nonumber \\ \label{eqIV-5}
&\sum\limits_{t=t}^{t+T-1} \!\!(P^{\mathrm{L, up}}_{t,i} - P_{t,i}^{\mathrm{L,do}})\Delta t, ~\forall i
\end{align}
Constraints (\ref{eqIV-4}) and (\ref{eqIV-5}) are used to ensure the satisfactions of (\ref{eqIII-25.5}) and (\ref{eqIII-27}) to balance the SoC of ESSs and satisfy the electricity consumption requirement of loads, respectively.

\subsection{Solution Methodology and Implementation}

In (\ref{eqIV-1}) and (\ref{eqIV-3}), there are many nonlinear terms in both objectives and constraints, which make the optimization problems difficult to solve. To tackle this, the following measures are adopted to convert these single-level problems into mixed integer second order cone programming (MISOCP) problems, thus the converted problem can be efficiently solved by available off-the-shelf solvers.

\emph{1) Linearization.}

Since the complementary relaxation constraints in (\ref{eqIII-33.1}) are nonlinear, thus the boolean variable $\bm{\kappa}$ is introduced to transform them into linear inequality with the form of:
\vspace{-3pt}
\begin{align}\label{eqIV-6}
&\textbf{0} \leq \bm{\mu} \leq \bm{M\kappa}  \\
&\textbf{0} \leq \bm{h(x)} \leq \bm{M(I-\kappa)}
\end{align}
where $M$ is a large enough positive number.

Besides the complementary relaxation constraints, variable substitution is utilized to realize the linearization of \emph{quadratic terms}, i.e., substituting $I_{t,ij}^2$ and $U_{t,i}^2$ with $l_{t,ij}$ and $v_{t,i}$, respectively. 
Doing this, Eqs. (\ref{eqIII-3})--(\ref{eqIII-5}), (\ref{eqIII-11}) and (\ref{eqIII-12}) can be transformed to linear constraints via variable substitution.

As for the nonlinear voltage deviation term in (\ref{eqIII-2.5}), auxiliary variable $Aux_{t,i}$ is introduced to express the extent of voltage deviation. The voltage deviation term can be thus linearized as:
\vspace{-8pt}
\begin{align}\label{eqIV-7}
&\mathcal{F}_{\mathrm{vd}} = \sum\limits_{t=1}^{N_T} \sum\limits_{i=1}^{N_N} {Aux_{t,i}}    \\
&Aux_{t,i} \geq v_{t,i}-\widetilde{U}_{\mathrm{ref}}^2; Aux_{t,i} \geq \widetilde{U}_{\mathrm{ref}}^2 - v_{t,i}; Aux_{t,i} \geq 0   \nonumber
\end{align}

Similarly, the tap adjusting cost of the OLTC $f_{\mathrm{sw}}$ can be linearized by introducing auxiliary variables $O_{t,ij}^{+}$ and $O_{t,ij}^{-}$, representing the positive/negative changes in the tap steps of the OLTC, respectively. In doing so, the cost calculation and operation constraints of the OLTC are rewritten as:
\vspace{-3pt}
\begin{align}\label{eqIV-8}
&f_{\mathrm{sw}} = \sum\limits_{ij\in \Omega_O}\sum\limits_{t=1}^{N_T}(O_{t,ij}^{+}+O_{t,ij}^{-});    \\
& \sum\nolimits_{t=1}^{N_T}(O_{t,ij}^{+}+O_{t,ij}^{-}) \leq \overline{\Delta}^{\mathrm{OLTC}}; \nonumber \\
& O_{t,ij}^{+}\geq 0; ~~O_{t,ij}^{-}\geq 0. \nonumber
\end{align}

After the variable substitution of $v_{t,i}$, constraint (\ref{eqIII-17.1}) can be re-expressed and transformed as follows. 
\vspace{-3pt}
\begin{align}\label{eqIV-9}
& v_{t,j} = v_{t,i}(1+O_{t,ij}\Delta o_{ij})^2, ~\forall t    \\
& O_{t,ij} = \sum\nolimits_{x=0}^{2\overline{O}_{ij}}(x-\overline{O}_{ij})\delta_{t,ij,k}, ~\delta_{t,ij,x}\in\{0,1\}   \\
& \sum\nolimits_{x=0}^{2\overline{O}_{ij}}\delta_{t,ij,x} = 1
\end{align}
where $\delta_{t,ij,x}$ is a binary variable to represent integer variable $O_{t,ij}$. Furthermore, on the basis of {\color{black}(\ref{eqIV-9}) and (45)}, we now have\vspace{-3pt}
\begin{align}\label{eqIV-10}
& v_{t,j} = \sum\nolimits_{x=0}^{2\overline{O}_{ij}}[(o_{ij,0}+(x-\overline{O}_{ij})\Delta o_{ij})^2 v_{t,ij,x}]    \\
& \underline{U}^2\delta_{t,ij,x} \leq v_{t,ij,x} \leq \overline{U}^2\delta_{t,ij,x}  \\
& \underline{U}^2(1-\delta_{t,ij,x}) \leq v_{t,j}-v_{t,ij,x} \leq \overline{U}^2(1-\delta_{t,ij,x})
\end{align}
where $v_{t,ij,x}$ is a auxiliary voltage variable to represent the nonlinear product of $v_{t,j}v_{t,ij,x}$.

\emph{2) Conic Relaxation.}

\emph{To quadratic cone constraints:} Using convex relaxation, Eq. (\ref{eqIII-6}) can be conversed to a quadratic cone constraint.
\vspace{-2pt}
\begin{align}\label{eqV-4}
\left\|\!\!\!
\begin{array}{c}
2P_{t,ij}  \\
2Q_{t,ij}  \\
l_{t,ij} - v_{t,i}
\end{array} \!\!\!\right\|_2 \leq l_{t,ij} + v_{t,i}, ~\forall t
\end{align}

\emph{To rotated cone constraints:} Similarly, Eqs. (\ref{eqIII-14}), (\ref{eqIII-16}) and (\ref{eqIII-10}) can be conversed to rotated quadratic cone constraints.
\vspace{-2pt}
\begin{equation}\label{eqV-5}
(P_{t,i}^{\mathrm{sop}})^2+(Q_{t,i}^{\mathrm{sop}})^2 \leq 2\frac{P_{t,i}^{\mathrm{sop,loss}}}{\sqrt{2}A_i^{\mathrm{sop}}}\frac{P_{t,i}^{\mathrm{sop,loss}}} {\sqrt{2}A_i^{\mathrm{sop}}} ~~~~~~~~~~~~
\end{equation}
\vspace{-5pt}
\begin{equation}\label{eqV-6}
(P_{t,i}^{alt})^2+(Q_{t,i}^{alt})^2 \leq 2\frac{S_i^{alt}}{\sqrt{2}}\frac{S_i^{alt}}{\sqrt{2}}, ~alt\in\{\mathrm{sop, p, w}\}
\vspace{-2pt}
\end{equation}

Through above transformation processes, i.e., linearizations and conic relaxation, the formulated non-linear programming problems (\ref{eqIV-1}) and (\ref{eqIV-3}) are transformed into MISOCP models, which can be efficiently solved using available optimization packages, such as Gurobi and CPLEX solvers.

Following the details mentioned above, the real-time control decisions of controllable resources in the DN and each MG can be determined using the proposed NTC method, which are implemented by \emph{Algorithm \ref{Algorithm1}}.

\begin{algorithm}[!h]
\caption{\color{black}{Implementation Algorithm of Proposed NTC}}
\label{Algorithm1}
\begin{algorithmic}[1]
\small
\State Initial parameters.
\State Formulate leader-followers transactive market between DNO and MGs according to (\ref{eqIII-32}).
\State Transform (\ref{eqIII-32}) to single level problem (\ref{eqIII-33}) via KKT conditions and duality.

\noindent Start to implement two-timescale NTC method in the following.

\noindent $~\diamond$ \underline{\emph{Slow-Timescale Pre-Scheduling Model}:}
\For {$t=1:N_{\mathrm{T}}$}
\State Collect $\&$ predict related information for $\{t, t+1, N_T\}$.
\State Formulate optimization problem (\ref{eqIV-1}).
\State Execute problem transformation on problem (\ref{eqIV-1}).
\State Solve, and deliver the results to \emph{real-time scheduling stage}.

\noindent $~\diamond$ \underline{\emph{ Real-Time Scheduling Model}:}
\While {$h\in[(t-1)\zeta, ~t\zeta]$}
\State Collect $\&$ short-term predict related information for $\textbf{\emph{T}}_{\mathrm{\textbf{hor}}}$.
\State Formulate optimization problem (\ref{eqIV-3}).
\State Execute problem transformation on problem (\ref{eqIV-3}).
\State Once the solution is obtained, implement only the actions associated with the current period $h$ to the physical system.
\State $h\leftarrow h+1$.
\EndWhile
\State $t\leftarrow t+1$.
\EndFor

\noindent $~\diamond$ \underline{\emph{Problem Transformations to MISOCP} (lines 18-21):}
\State Linearize objective function through variable substitution and absolute term elimination.
\State Transform constraints into linear or conic constraints.
\State MISOCP problem is thus re-formulated from original problem.

\State \textbf{Output:} real-time scheduling strategies for OLTC, SOPs, ESSs and DR activities.
\label{code:recentEnd}
\end{algorithmic}
\end{algorithm}

\section{Simulation Studies} \label{sec5}

\subsection{Simulation Setups}\label{sec.5.1}

\begin{figure}[!t]
  \centering
  \includegraphics[width=3.4in, height=1.4in]{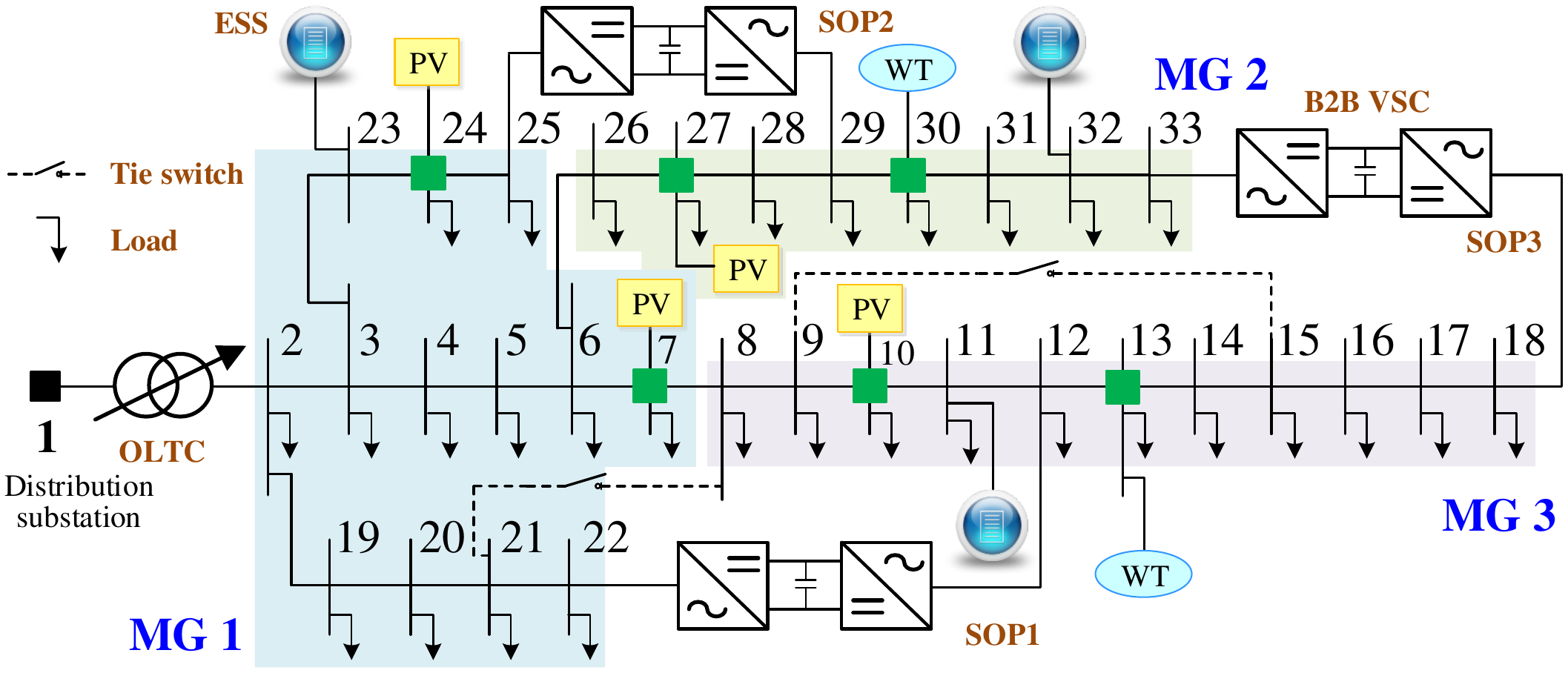} \\
  \caption{Topology diagram of the modified IEEE 33-bus test system with three microgrids and three SOPs.}
  \label{figfSecV-1}
\vspace{-10pt}
\end{figure}

\begin{figure}[!t]
  \centering
  \includegraphics[width=3.0in]{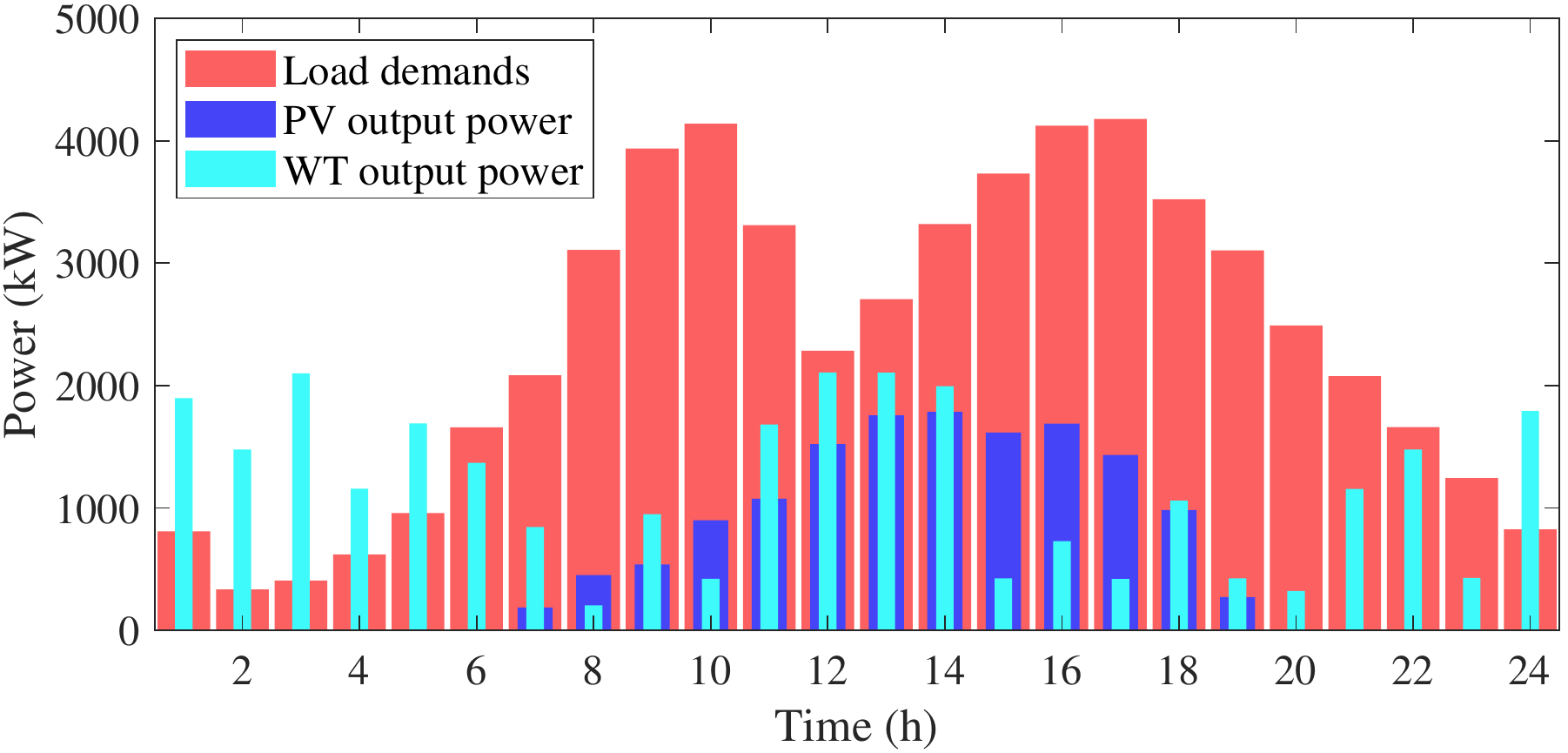} \\
  \caption{Original predicted data for load demands, PV and WT output profiles.}
  \label{figfSecV-2}
\vspace{-10pt}
\end{figure}

\emph{1) Network Model:} A modified IEEE 33-bus test system with three microgrids is used to validate the effectiveness of our proposed NTC method, as shown in Fig. \ref{figfSecV-1}. An OLTC, four PV, two WT generators, three ESSs and three SOPs are additionally integrated, and their corresponding placements and parameters are shown in Table I. The other parameters are set as the standard IEEE 33-bus system.
The time-series of system general loads, WT and PV active power are depicted in Fig. \ref{figfSecV-2}, based on which the utilized scenarios can be generated via (\ref{eqA35}).

\begin{table}[!t]
\scriptsize
\newcommand{\tabincell}[2]{\begin{tabular}{@{}#1@{}}#2\end{tabular}}
\caption{Modifications in IEEE 33-bus test system}
\label{tabVI-1}
\centering
\begin{tabular}{l|l|l|l}
\hline
\!\textbf{Entity}\! & \textbf{Device} & \textbf{Location} & \textbf{Parameter}  \\
\hline
\multirow{3}{*}{MG 1} & PV unit & \tabincell{c}{Buses $7$, $24$} & \tabincell{c}{$650$kW, $400$kW}  \\
\cline{2-4}
 & Energy storage  & Bus 6 & \tabincell{l}{1.0MWh, 0.2MW, 0.95}  \\
\hline
\multirow{3}{*}{MG 2} & PV unit & \tabincell{c}{Bus $27$} & \tabincell{c}{$500$kW}  \\
\cline{2-4}
 & WT unit & Bus $30$  &  \tabincell{c}{$1600$kW} \\
\cline{2-4}
 & Energy storage  & Bus 32 & \tabincell{l}{1.0MWh, 0.2MW, 0.95}  \\
\hline
\multirow{3}{*}{MG 3} & PV unit & \tabincell{c}{Bus $10$} & \tabincell{c}{$650$kW}  \\
\cline{2-4}
 & WT unit & Bus $13$  &  \tabincell{c}{$1300$kW} \\
\cline{2-4}
 & Energy storage  & Bus 16 & \tabincell{l}{1.0MWh, 0.2MW, 0.95}  \\
\hline
\multirow{2}{*}{DN} & SOPs & \!\tabincell{l}{Buses 12-22, 25-29, \\18-33}\! & \tabincell{l}{Capacity: $1.0$MVA; \\Control mode: $PQ$-$V_{\mathrm{dc}}Q$}  \\
\cline{2-4}
 & OLTC & Buses 1-2  &  \tabincell{c}{$\pm5 \times1\%$ (include $0\times1\%$)} \\
\hline
\end{tabular}
\vspace{-10pt}
\end{table}

\emph{2) Scheduling Parameter Settings:} The time intervals of two timescales are 1 hour and 1/4 hour, respectively. The capacity of the installed SOPs is set as $1.0$ MVA. The percentage of demands that participates the DR program is assumed to be 20\% \cite{MultiFlexibility2020Marco,yang2019interactive}. The prices that the DN buys electricity from the HV system are set as Ref. \cite{li2018energypricing}, and the selling price of the DN is set as {\color{black}0.0584 $\$$/kWh}.
All the installed renewable generators are operated at a unit power factor without considering the localized reactive power support of renewables \cite{yang2020flexibility}.
Other scheduling parameter settings are given in Table \ref{tab2}.

\begin{table}[!t]
\scriptsize
\newcommand{\tabincell}[2]{\begin{tabular}{@{}#1@{}}#2\end{tabular}}
\caption{Parameters settings}
\label{tab2}
\centering
\begin{tabular}[!t]{|cc|cc|cc|}
  \hline
  \!\!{\textbf{Parameter}}\!\! & \textbf{Value} & \!\!{\textbf{Parameter}}\!\! & \textbf{Value} & \!\!{\textbf{Parameter}}\!\! & \textbf{Value}  \\
  \hline
  $\Delta t$ & 1 hour & $\Delta \kappa$ & $1/4$ hour & \tabincell{c}{$\overline{\Delta}^{\mathrm{OLTC}}$} & \!\!{4 times$/$day}\!\!  \\
  \tabincell{c}{$\mathbbm{a}_{\mathrm{o}}$} & 0.833 & $\mathbbm{b}_{\mathrm{v}}$ & 0.167\cite{li2017coordinated} & $\overline{U}$ & \!{1.05 p.u.}\! \\
  \tabincell{c}{$\mathscr{C}_{\mathrm{loss}}$} & \!\!{Same as $X_{t}$}\!\! & $\mathscr{C}_{\mathrm{deg}}$ & \!\!{\color{black}2.736 $\$$/MWh}\!\! & \tabincell{c}{$\underline{U}$} & \!{0.95 p.u.}\!  \\
  $o_{ij,0}$  & $1.0$  & $\Delta o_{ij}$ & $1\%$ & $\widetilde{U}_{\mathrm{ref}}$ & 1.0 p.u.   \\
  $\lambda_t^{\mathrm{min}}$  & $0.8X_{t}$  & $\lambda_t^{\mathrm{max}}$ & $1.2X_{t}$ & $T$ & 3 hours   \\
  \hline
\end{tabular}
\vspace{-10pt}
\end{table}

\emph{3) Simulation Environment:} All the simulations are conducted in the MATLAB 2018b in a 64-bit Windows environment with Gurobi $9.0$ solver and YALMIP toolbox, on a PC with Core i7-8700 CPU @3.2 GHz processor and 8 GB RAM.

\subsection{Results and Analysis}

This part aims to verify the reasonability of the proposed NTC method from the respects of transactive control results, statistical performance and voltage magnitude profile.

\begin{figure}[!t]
  \centering
  \includegraphics[width=3.5in]{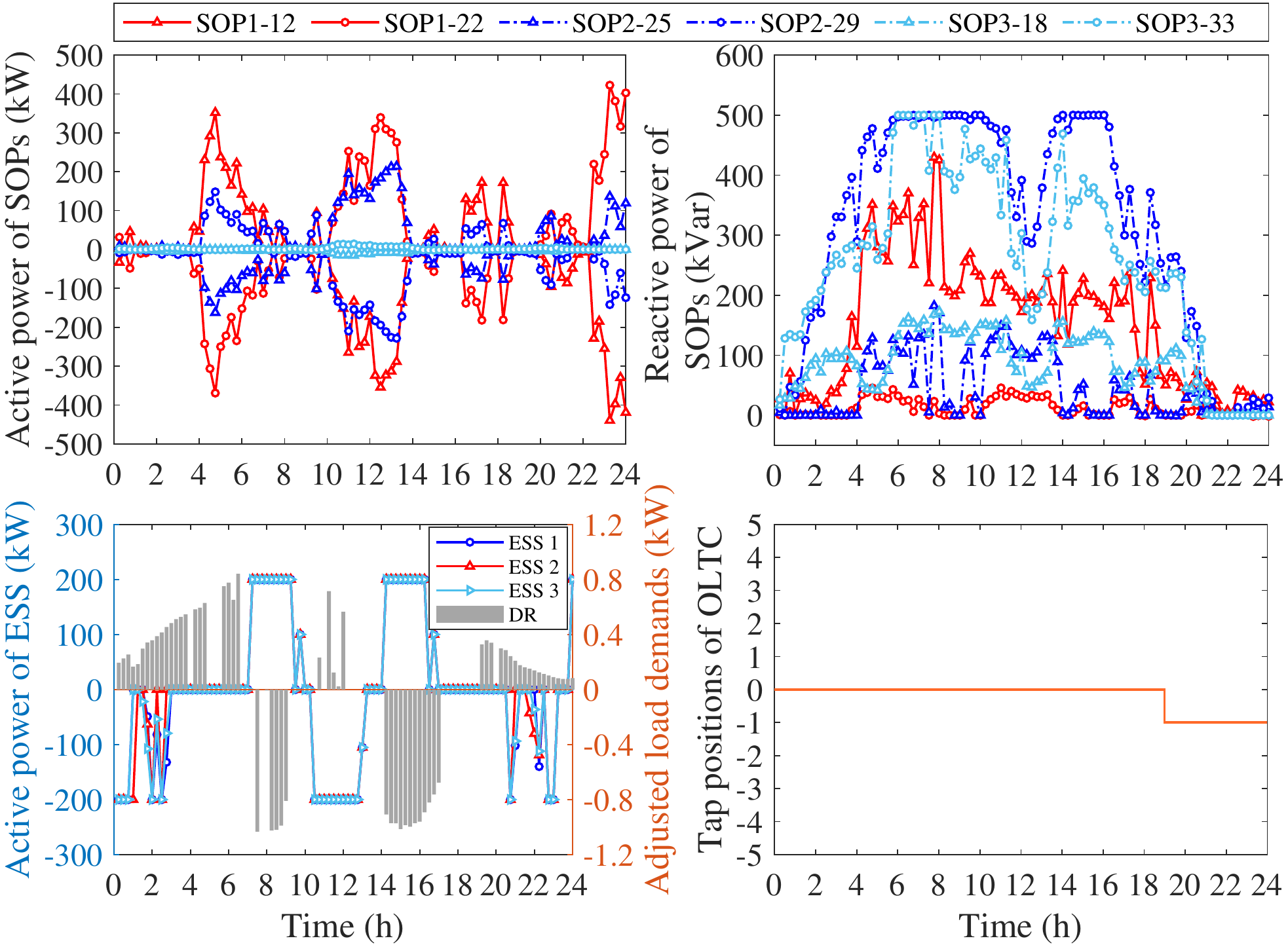}  \\
  \caption{Final transactive control results: (a)-(b) active and reactive power of the three SOPs; (c) control strategies of ESSs and DR activities; (d) optimal tap positions of the OLTC.}
  \label{figfSecV-3}
\vspace{-10pt}
\end{figure}

\begin{figure}[!t]
  \centering
  \includegraphics[width=3.0in]{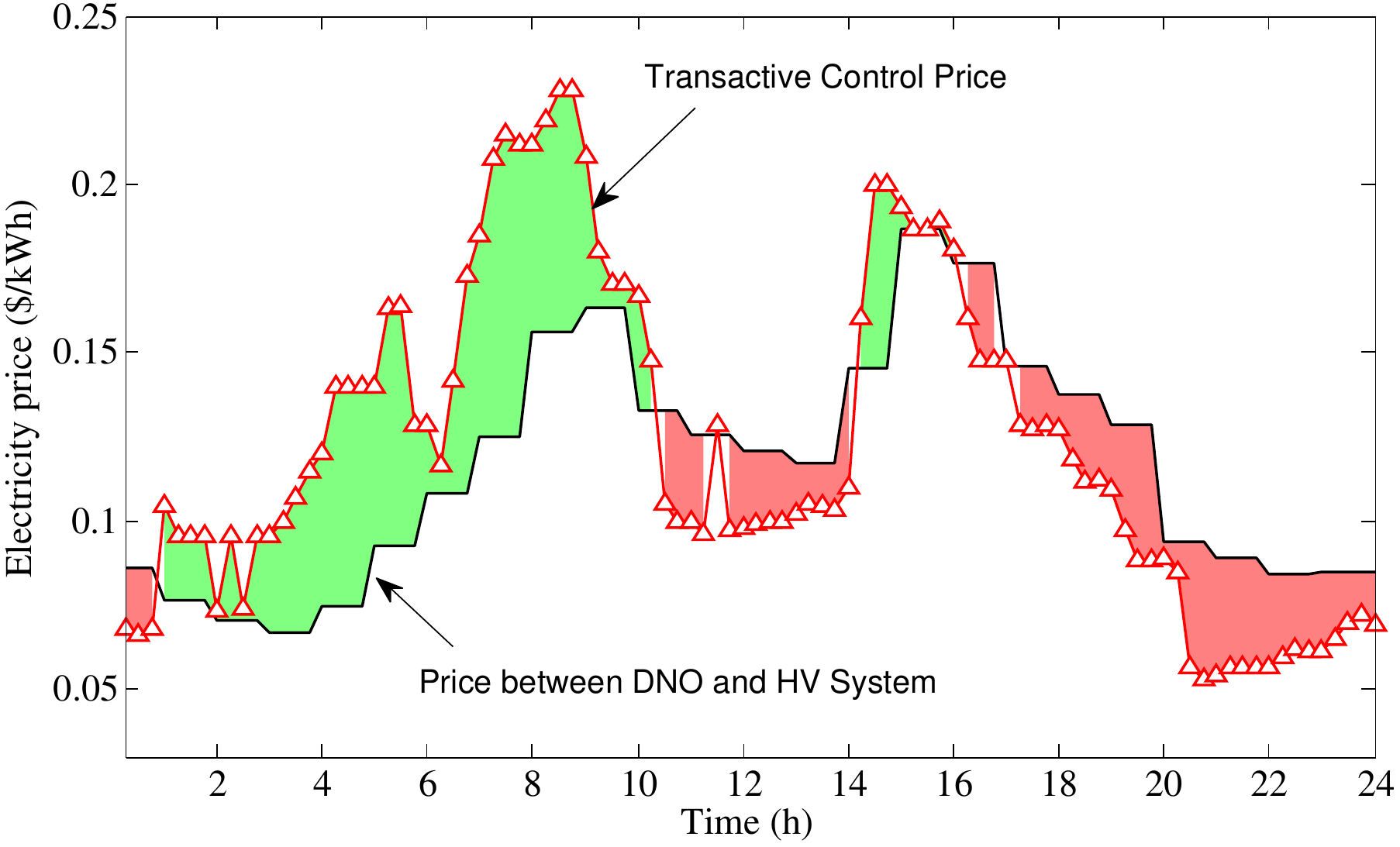}  \\
  \caption{Transactive energy price curve between DNO and MGs.}
  \label{figfSecV-4}
\vspace{-10pt}
\end{figure}

\begin{figure}[!t]
  \centering
  \includegraphics[width=3.0in]{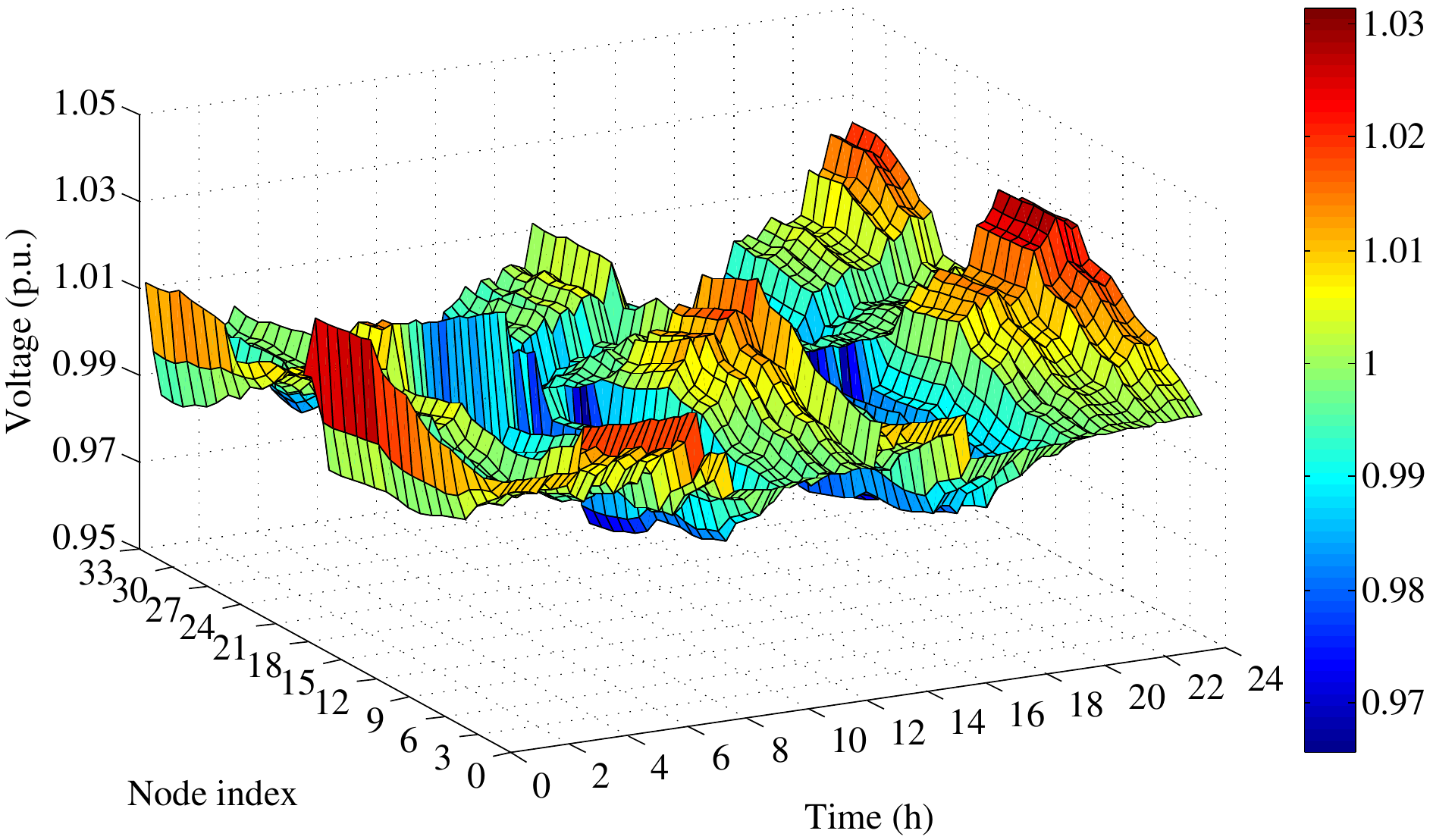}  \\
  \caption{Voltage profiles of all 33 nodes over an operational day.}
  \label{figfSecV-5}
\vspace{-10pt}
\end{figure}

\emph{1) Transactive control results of controllable resources.} By implementing the NTC method on the modified IEEE 33-bus test system with previous setups, the obtained transactive control results concerning the active/reactive power of SOPs, charding/discharging power of ESSs, DR activities and tap positions of OLTC are shown in Fig. \ref{figfSecV-3}. The transactive energy price profile between the DNO and MGs is given in Fig. \ref{figfSecV-4}.

\emph{2) Performance analysis.} Numerical performances of the test 33-bus system under different optimization situations are listed in Table \ref{tab3}.
From the results in Figs. \ref{figfSecV-3}-\ref{figfSecV-4} and Table \ref{tab3}, we can see that, the ESSs discharge and flexible loads reduce consumption during the renewable power shortage periods such as 7:00-10:00 and 14:00-20:00, and quite the contrary during the renewable power is abundant (i.e., 00:00-5:00 and 11:00-14:00), thus contributing to improving the local supply-demand balance.
Moreover, SOPs cooperate with the OLTC to timely respond the voltage volatility caused by RESs, which can effectively lower the security risks of the DN operation. With accurately power flow regulating ability, SOPs help to greatly reduce line losses and voltage deviations compared to the unscheduling mode.

\begin{table}[!t]
\scriptsize
\newcommand{\tabincell}[2]{\begin{tabular}{@{}#1@{}}#2\end{tabular}}
\caption{Performances data of the MMG-based distribution network}
\label{tab3}
\centering
\begin{tabular}{c|c|c|c|c}
\hline
\multirow{3}{*}{System performance} & \!\!\multirow{3}{*}{\tabincell{c}{Unscheduling\\mode$^1$}}\!\! & \multicolumn{2}{c|}{Proposed NTC method} & \multirow{3}{*}{\tabincell{c}{Ideal\\mode$^2$}}   \\
\cline{3-4}
 & & \!\tabincell{c}{Pre-\\Scheduling}\!  &  \!\tabincell{c}{Real-time \\Scheduling}\! &  \\
\hline
\tabincell{c}{Line power losses (kWh)} & 1736.04 & 756.86  & 749.82 & 718.62  \\
\tabincell{c}{SOP power losses (kWh)} & 0 & 502.50  & 509.97 & 518.39  \\
\!\tabincell{c}{Voltage deviation (p.u.$^2$)\\(calculated by $\mathcal{F}_{\mathrm{vd}}$ in (\ref{eqIII-2.5}))}\! & 22.86 & 18.00  & 13.57 & 12.96  \\
\hline
\tabincell{c}{Profit of DNO ($\$$)} & 430.28 & 635.176  & 753.907 & 790.53  \\
\tabincell{c}{Total cost of MGs ($\$$)} & 5536.96 & 5324.08  & 5045.86 & 4825.85  \\
\hline
\end{tabular}
\begin{tablenotes}
  \vspace{2pt}
  \item[1] $^1$For the unscheduling mode, it simulates the passive operation of the test system, where the optimal control of the flexible resources belonging to the DNO and MGs is not implemented (a reference case). Herein, the voltage limits are removed to avoid an unacceptably huge amount of load shedding and renewable curtailment \cite{yang2020flexibility}.
  \vspace{4pt}
  \item[1] $^2$For the ideal mode, it is assumed that all the needed information are known to DNO in advance with perfect predictions.
\end{tablenotes}
\vspace{-10pt}
\end{table}

Also, Table \ref{tab3} further shows that the performances under real-time scheduling are obviously better than the ones under pre-scheduling. The reason lies in the fact that the final control strategies are obtained via a dynamic updated manner in smaller scheduling granularity based on the latest updated information, which has proved to be efficient enough to against system uncertainties.

\emph{3) Voltage profile.} System voltage magnitude is an important indicator to represent the operational security of DNs. Thus, the voltage profiles of all the 33 nodes over real-time scheduling are shown in Fig. \ref{figfSecV-5}. Obviously, the proposed NTC can eliminate all voltage violations by adjusting system power flow via SOPs, thus strictly maintaining the voltage magnitudes within safety range (\emph{i.e.}, [0.95, 1.05] p.u.). As such, the security of the DN operation can be ensured. Moreover, compared to the unscheduling mode, the voltage deviation performed by the NTC method is reduced by {\color{black}$43.30\%$}, which could mitigate the adverse effects of the uncertainties on system security.

According to above analyses, the proposed NTC method employs the SOPs as an active measure in DN-level for supplementing the transactive market between DNO and MMGs, which is proved to be of great ability to eliminate voltage violations and decrease power losses, thus effectively improving the operational security and economy of DNs simultaneously.

\subsection{Sensitive Analysis Towards Uncertainties}

\begin{figure}[!t]
  \centering
  \includegraphics[width=3.3in]{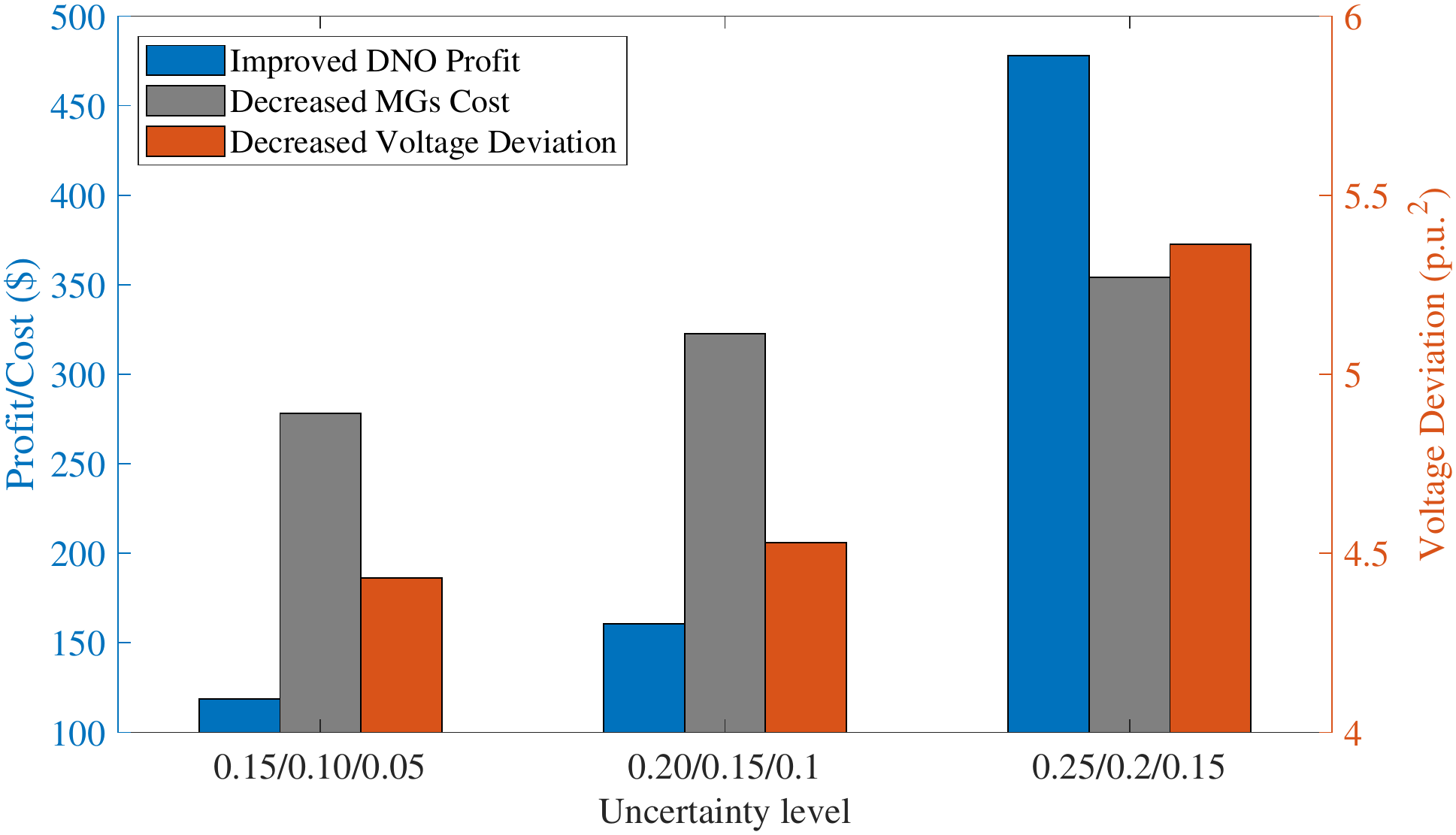}  \\
  \caption{The improved performances of fast-timescale stage. (Note: the uncertainty level $a/b/c$ means that uncertainties $u_{i,r}$ of WT, PV, and load demands are set as $a$, $b$ and $c$, respectively.)}
  \label{figfSecV-6}
\vspace{-10pt}
\end{figure}

In reality, the predictions are always not perfect due to the inherent intermittency and volatility of RESs \cite{ding2018data,yang2019event}. Therefore, in order to verify the robustness of our proposed NTC method for working against uncertainties, a sensitive analysis under various uncertain settings is performed in this part.
Consequently, Table \ref{tab4} lists the test results under different uncertainty levels of prediction errors, and Fig. \ref{figfSecV-6} shows the improved performances of real-time scheduling stage with comparison of slow-timescale stage in three test scenarios.

\begin{table}[!t]
\scriptsize
\newcommand{\tabincell}[2]{\begin{tabular}{@{}#1@{}}#2\end{tabular}}
\caption{Performance comparison of NTC method under different uncertainty levels}
\label{tab4}
\centering
\begin{tabular}{ll|cccc}
\hline
\!\!\multirow{3}{*}{\tabincell{c}{Unc.\\Scena.}}\!\!\! & \multirow{3}{*}{Performances} & \!\!\multirow{3}{*}{\tabincell{c}{Deterministic\\mode$^1$}}\!\! & \multicolumn{2}{c}{Proposed NTC method} & \multirow{3}{*}{\tabincell{c}{Ideal\\mode}}   \\
\cline{4-5}
 & & & \!\tabincell{c}{Pre-\\Scheduling}\!  &  \!\tabincell{c}{Real-time \\Scheduling}\! &  \\
\hline
\multirow{4}{*}{1} & \!\!\tabincell{c}{Power losses (kWh)}\!\! & 919.63 & 756.86 & 749.82 & 718.62   \\
 & \tabincell{c}{Profit of DNO ($\$$)} & 554.18 & 635.176 & 753.907  & 790.53   \\
 & \!\!\tabincell{c}{Tot. cost of MGs ($\$$)}\!\! & 5425.31 & 5324.08 & 5045.86  & 4825.85   \\
 & \tabincell{c}{Volt. dev. (p.u.$^2$)} & 20.24 & 18.00 & 13.57  & 12.96   \\
\hline
\multirow{4}{*}{2} & \!\!\tabincell{c}{Power losses (kWh)}\!\! & 927.38 & 765.08 & 754.37  & 709.10   \\
 & \tabincell{c}{Profit of DNO ($\$$)} & 556.90 & 589.608 & 750.131  & 801.81   \\
 & \!\!\tabincell{c}{Tot. cost of MGs ($\$$)}\!\! & 5422.32 & 5372.43 & 5049.79  & 4746.80   \\
 & \tabincell{c}{Volt. dev. (p.u.$^2$)} & 20.37 & 18.17 & 13.64  & 12.80   \\
\hline
\multirow{4}{*}{3} & \!\!\tabincell{c}{Power losses (kWh)}\!\! & 933.56 & 746.56 & 722.90  & 663.91   \\
 & \tabincell{c}{Profit of DNO ($\$$)} & 452.81 & 456.335 & 934.126  & 1024.64   \\
 & \!\!\tabincell{c}{Tot. cost of MGs ($\$$)}\!\! & 5388.19 & 5365.17 & 5011.11  & 4612.23   \\
 & \tabincell{c}{Volt. dev. (p.u.$^2$)} & 20.39 & 18.68 & 13.32  & 12.25   \\
\hline
\end{tabular}
\begin{tablenotes}
  \item[1] {\color{black}$^1$ A deterministic NTC formulation using predicted information for decision-making without considering uncertainties.}
\end{tablenotes}
\vspace{-5pt}
\end{table}

\begin{table*}[!t]
\scriptsize
\newcommand{\tabincell}[2]{\begin{tabular}{@{}#1@{}}#2\end{tabular}}
\caption{Performance Comparison Among the Five Cases}
\label{tab5}
\centering
\begin{tabular}{c|c|c|c|c|c}
\hline
\tabincell{c}{System performance} & \tabincell{c}{Unscheduling mode} & ~~~~~\tabincell{c}{Case I}~~~~~ & ~~~~~\tabincell{c}{Case II}~~~~~ & ~~~~~\tabincell{c}{Case III}~~~~~ & \tabincell{c}{Case IV (proposed NTC method)}  \\
\hline
\tabincell{c}{Line power losses (kWh)} & 2056.486 & 908.736 & 819.074 & 1942.462 & 769.327    \\
\tabincell{c}{SOP power losses (kWh)} & 0 & 632.962 & 545.477 & 0 & 526.263    \\
\tabincell{c}{Voltage deviation (p.u.$^2$)} & 23.460 & 89.471 & 15.860 & 23.381 & 14.228    \\
\hline
\tabincell{c}{Profit of DNO ($\$$)} & 418.821 & 830.446 & 1353.810 & 802.722 & 919.325    \\
\tabincell{c}{Total cost of MGs ($\$$)} & 5219.757 & 4973.600 & 5421.541 & 5093.116 & 5018.813    \\
\hline
\end{tabular}
\vspace{-10pt}
\end{table*}

From the results we can know that, with the increasing levels of prediction errors, the various performances of unscheduling model consistently become worse.
Also, Fig. \ref{figfSecV-6} shows that the improved performances in terms of DNO profit, total MGs cost and voltage deviation are more visible with the ascending levels of prediction errors, but the gap from the ideal mode is still constantly widening.

These facts indicate that the prediction uncertainties concerning renewables and consumptions will have adverse effects on the economic and secure performances of the DN system, and this kind of adverse effects increase with the prediction error increases.
Furthermore, compared to unscheduling and single-layer models, the robustness of the proposed two-timescale NTC method is verified as it can indeed mitigate the adverse effects of the uncertainties through dynamic double-layer strategies updating, and it shows more potential in cost savings and voltage regulation as uncertainties grow larger.

\subsection{Case studies}

In this subsection, we tend to validate the effectiveness of the proposed NTC method and examine the effect of each salient feature through case studies. For this purpose, four typical cases are further conducted and compared, which are defined as follows.
We randomly choose a set of scenarios under uncertainty level of {\color{black}$0.25/0.2/0.15$} for test.

\begin{figure}[!t]
  \centering
  \includegraphics[width=3.0in]{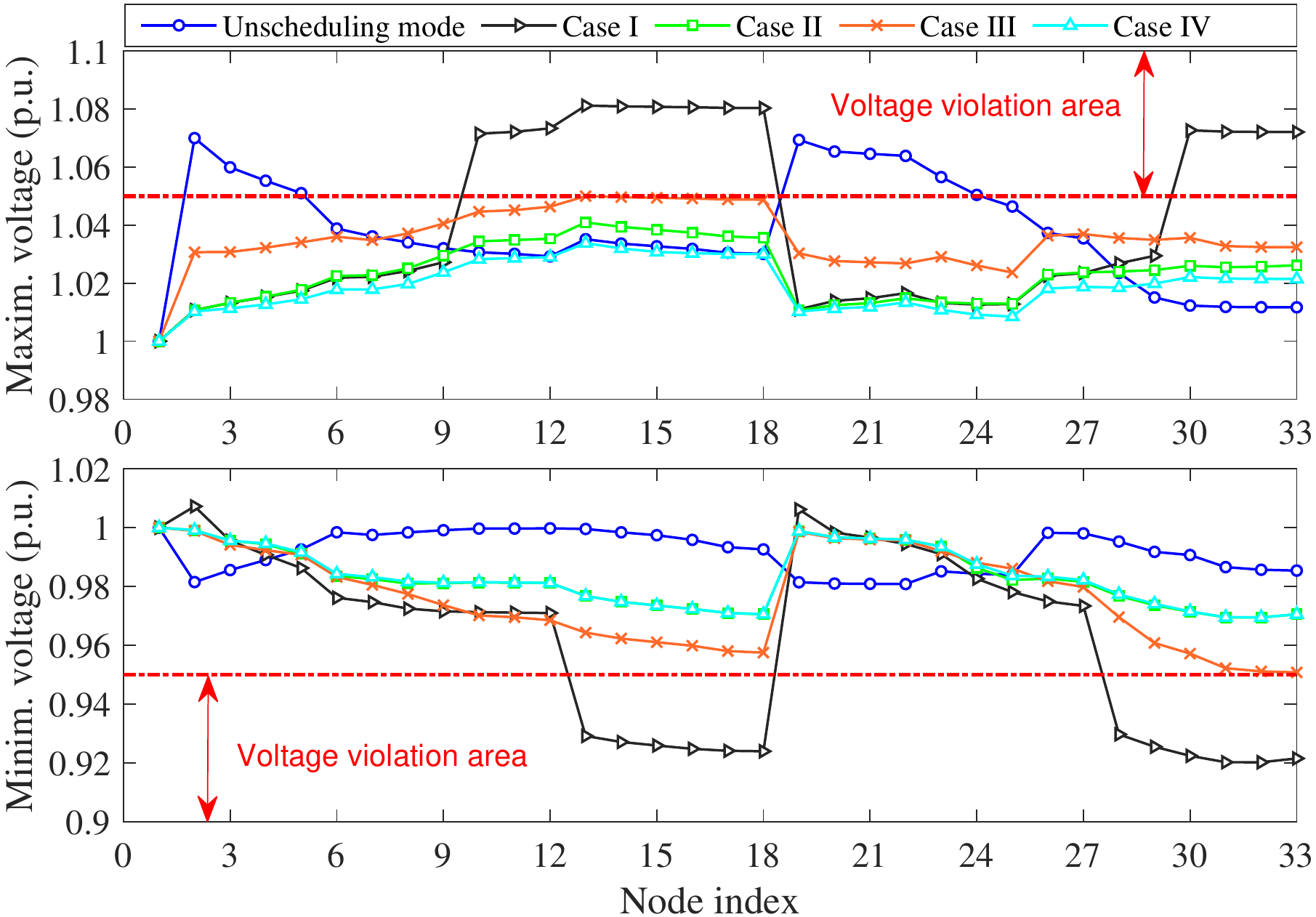} \\
  \caption{\color{black}Extreme voltage magnitude comparisons of the 5 cases. (a) Maximum voltage profiles; (b) Minimum voltage profiles in each node.}
  \label{figfSecV-7}
\vspace{-10pt}
\end{figure}

\begin{figure}[!t]
  \centering
  \includegraphics[width=3.1in]{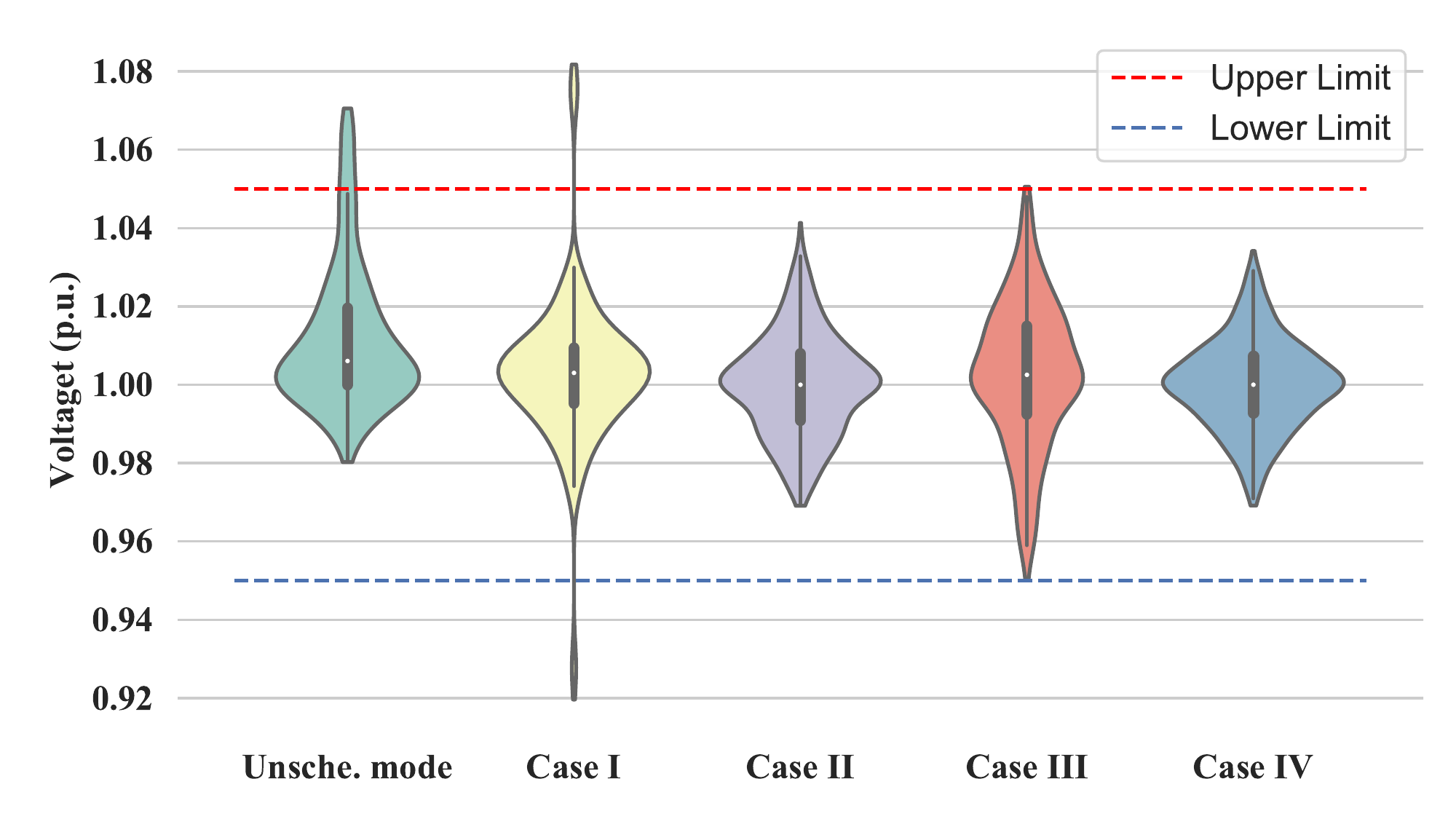} \\
  \caption{\color{black}Voltage magnitude distributions in Violin plot for 5 cases.}
  \label{figfSecV-8}
\vspace{-10pt}
\end{figure}


\emph{Case I:} Economy-oriented transactive control. In this case, only economy-concerned factors are considered in the transactive control framework, and the power flow, line power losses and voltage deviation are not included.

\emph{Case II:} Network-constrained non-market-based control. In this case, the trading price between the DNO and MGs is the HV grid price multiplied by a number greater than 1 (i.e., $1.2$). That is, the transactive energy market is not considered.

\emph{Case III:} Network-constrained transactive control without considering the active measure of controlling SOPs.

\emph{Case IV:} The proposed NTC method in previous sections.

With above definitions, the comparison results of five cases (including unscheduling mode) are listed in Table \ref{tab5}. Their extreme voltage magnitude comparisons and voltage magnitude distributions in violin plot are shown in Fig. \ref{figfSecV-7} and Fig. \ref{figfSecV-8}, respectively.

From the results, it can be seen that our proposed NTC method performs better than other four cases in both economy and security aspects. Case I is an economy-oriented method, and thus results in a lower MGs cost due to the fully utilization of ESS and DR without considering the violations of network constraints. We further conduct power lower calculation for it using its obtained strategies, the results show that it has larger line losses and voltage deviation in comparison with Case IV, and raises voltage violations as depicted in Fig. \ref{figfSecV-8}.
In Case II, drooping transactive energy market results in uneven distribution of benefits between DNO and MGs. In this case, the ESS and DR action less, so the SOPs in DN-level make more efforts for power flow regulation, thus leading to larger SOP power losses.
In Case III, after removing SOPs, the line losses and voltage deviation are increased by {\color{black}$152.49\%$ and $64.33\%$} respectively compared to Case IV, which precisely verifies the benefits of SOPs in reshaping system power flow to deal with security issues of DNs.

In summary, we can conclude from above analyses that the proposed NTC method is able to actively optimize the secure operation of the MMG-based DN system besides providing a platform to coordinate the DNO and MGs with respectation of their interests, through innovatively integrate a transactive energy market with the novel power-electronics device SOP-based power flow regulation technique.
Although a fraction of economy is sacrificed, the security of entire system operation can be thus absolutely guaranteed.

\section{Conclusion} \label{sec6}

In this paper, we propose a network-constrained transactive control framework for an MMG-based distribution network considering uncertainties. Different from previous TC studies, this framework can not only address the economic issues of transactive market between the DNO and the MGs, but also can adjust DN operation by optimally regulating the OLTC and the novel power-electronics device (i.e., SOPs). In this way, economic and technical issues are thus addressed in a holistic manner.
In particular, we innovatively integrate a collaborative optimization mechanism with the OPF technique.
A dynamic two-timescale model is formulated to minimize operational cost, improve voltage profile and against adverse effects of uncertainties.
Moreover, several model transformation and relaxation techniques are applied to avoid iterative solving process.
Case studies on {\color{black}a modified IEEE 33-bus distribution feeder with three MGs} demonstrate the effectiveness of our proposed NTC method and algorithm.

\bibliography{xd_yang_TII}

\begin{thebibliography}{10}
\providecommand{\url}[1]{#1}
\csname url@samestyle\endcsname
\providecommand{\newblock}{\relax}
\providecommand{\bibinfo}[2]{#2}
\providecommand{\BIBentrySTDinterwordspacing}{\spaceskip=0pt\relax}
\providecommand{\BIBentryALTinterwordstretchfactor}{4}
\providecommand{\BIBentryALTinterwordspacing}{\spaceskip=\fontdimen2\font plus
\BIBentryALTinterwordstretchfactor\fontdimen3\font minus
  \fontdimen4\font\relax}
\providecommand{\BIBforeignlanguage}[2]{{%
\expandafter\ifx\csname l@#1\endcsname\relax
\typeout{** WARNING: IEEEtran.bst: No hyphenation pattern has been}%
\typeout{** loaded for the language `#1'. Using the pattern for}%
\typeout{** the default language instead.}%
\else
\language=\csname l@#1\endcsname
\fi
#2}}
\providecommand{\BIBdecl}{\relax}
\BIBdecl

\bibitem{feng2020coalitional}
C.~{Feng}, F.~{Wen}, S.~{You}, Z.~{Li}, F.~{Shahnia}, and M.~{Shahidehpour},
  ``Coalitional game-based transactive energy management in local energy
  communities,'' \emph{IEEE Trans. Power Syst.}, vol.~35, no.~3, pp.
  1729--1740, 2020.

\bibitem{shuai2020online}
H.~{Shuai} and H.~{He}, ``Online scheduling of a residential microgrid via
  monte-carlo tree search and a learned model,'' \emph{IEEE Trans. Smart Grid},
  2020, to be published. DOI: \textcolor[rgb]{0,0,1}{10.1109/TSG.2020.3035127}.

\bibitem{xu2020MultiEnergy}
D.~{Xu}, B.~{Zhou}, N.~{Liu}, Q.~{Wu}, N.~{Voropai}, C.~{Li}, and
  E.~{Barakhtenko}, ``Peer-to-peer multi-energy and communication resource
  trading for interconnected microgrids,'' \emph{IEEE Trans. Ind. Informat.},
  2020, to be published. DOI: \textcolor[rgb]{0,0,1}{10.1109/TII.2020.3000906}.

\bibitem{daneshvar2020twostage}
M.~{Daneshvar}, B.~{Mohammadi-Ivatloo}, K.~{Zare}, and S.~{Asadi}, ``Two-stage
  robust stochastic model scheduling for transactive energy based renewable
  microgrids,'' \emph{IEEE Trans. Ind. Informat.}, vol.~16, no.~11, pp.
  6857--6867, 2020.

\bibitem{yang2019interactive}
X.~Yang, H.~He, Y.~Zhang, Y.~Chen, and G.~Weng, ``Interactive energy management
  for enhancing power balances in multi-microgrids,'' \emph{IEEE Trans. Smart
  Grid}, vol.~10, no.~6, pp. 6055--6069, 2019.

\bibitem{lezama2019local}
F.~{Lezama}, J.~{Soares}, P.~{Hernandez-Leal}, M.~{Kaisers}, T.~{Pinto}, and
  Z.~{Vale}, ``Local energy markets: Paving the path toward fully transactive
  energy systems,'' \emph{IEEE Trans. Power Syst.}, vol.~34, no.~5, pp.
  4081--4088, 2019.

\bibitem{cui2020efficient}
S.~{Cui}, Y.~W. {Wang}, Y.~{Shi}, and J.~W. {Xiao}, ``An efficient peer-to-peer
  energy-sharing framework for numerous community prosumers,'' \emph{IEEE
  Trans. Ind. Informat.}, vol.~16, no.~12, pp. 7402--7412, 2020.

\bibitem{zhou2020flexible}
Q.~{Zhou}, M.~{Shahidehpour}, A.~{Alabdulwahab}, and A.~{Abusorrah}, ``Flexible
  division and unification control strategies for resilience enhancement in
  networked microgrids,'' \emph{IEEE Trans. Power Syst.}, vol.~35, no.~1, pp.
  474--486, 2020.

\bibitem{farzin2016enhancing}
H.~{Farzin}, M.~{Fotuhi-Firuzabad}, and M.~{Moeini-Aghtaie}, ``Enhancing power
  system resilience through hierarchical outage management in
  multi-microgrids,'' \emph{IEEE Trans. Smart Grid}, vol.~7, no.~6, pp.
  2869--2879, 2016.

\bibitem{zhao2018energy}
B.~Zhao, X.~Wang, D.~Lin, M.~M. Calvin, J.~C. Morgan, R.~Qin, and C.~Wang,
  ``Energy management of multiple microgrids based on a system of systems
  architecture,'' \emph{IEEE Trans. Power Syst.}, vol.~33, no.~6, pp.
  6410--6421, 2018.

\bibitem{jadhav2019novel}
A.~M. Jadhav, N.~R. Patne, and J.~M. Guerrero, ``A novel approach to
  neighborhood fair energy trading in a distribution network of multiple
  microgrid clusters,'' \emph{IEEE Trans. Ind. Electron.}, vol.~66, no.~2, pp.
  1520--1531, 2019.

\bibitem{anoh2020energy}
K.~{Anoh}, S.~{Maharjan}, A.~{Ikpehai}, Y.~{Zhang}, and B.~{Adebisi}, ``Energy
  peer-to-peer trading in virtual microgrids in smart grids: A game-theoretic
  approach,'' \emph{IEEE Trans. Smart Grid}, vol.~11, no.~2, pp. 1264--1275,
  2020.

\bibitem{zhao2020distributed}
Z.~{Zhao}, J.~{Guo}, C.~S. {Lai}, H.~{Xiao}, K.~{Zhou}, and L.~L. {Lai},
  ``Distributed model predictive control strategy for islands multi-microgrids
  based on non-cooperative game,'' \emph{IEEE Trans. Ind. Informat.}, 2020, to
  be published. DOI: \textcolor[rgb]{0,0,1}{10.1109/TII.2020.3013102}.

\bibitem{esfahani2019multiagent}
M.~M. Esfahani, A.~Hariri, and O.~A. Mohammed, ``A multiagent-based
  game-theoretic and optimization approach for market operation of
  multimicrogrid systems,'' \emph{IEEE Trans. Ind. Informat.}, vol.~15, no.~1,
  pp. 280--292, 2019.

\bibitem{park2016contribution}
S.~Park, J.~Lee, S.~Bae, G.~Hwang, and J.~K. Choi, ``Contribution-based
  energy-trading mechanism in microgrids for future smart grid: A game
  theoretic approach,'' \emph{IEEE Trans. Ind. Electron.}, vol.~63, no.~7, pp.
  4255--4265, 2016.

\bibitem{jin2020local}
X.~Jin, Q.~Wu, and H.~Jia, ``Local flexibility markets: Literature review on
  concepts, models and clearing methods,'' \emph{Appl. Energy}, vol. 261, p.
  114387, 2020.

\bibitem{wang2020reconfigurable}
Y.~{Wang}, Z.~{Huang}, M.~{Shahidehpour}, L.~L. {Lai}, Z.~{Wang}, and Q.~{Zhu},
  ``Reconfigurable distribution network for managing transactive energy in a
  multi-microgrid system,'' \emph{IEEE Trans. Smart Grid}, vol.~11, no.~2, pp.
  1286--1295, 2020.

\bibitem{ji2019robust}
H.~Ji, C.~Wang, P.~Li, F.~Ding, and J.~Wu, ``Robust operation of soft open
  points in active distribution networks with high penetration of photovoltaic
  integration,'' \emph{IEEE Trans. Sustain. Energy}, vol.~10, no.~1, pp.
  280--289, 2019.

\bibitem{ding2018data}
T.~Ding, Q.~Yang, Y.~Yang, C.~Li, Z.~Bie, and F.~Blaabjerg, ``A data-driven
  stochastic reactive power optimization considering uncertainties in active
  distribution networks and decomposition method,'' \emph{IEEE Trans. Smart
  Grid}, vol.~9, no.~5, pp. 4994--5004, 2018.

\bibitem{li2019optimaloperation}
P.~{Li}, H.~{Ji}, C.~{Wang}, J.~{Zhao}, G.~{Song}, F.~{Ding}, and J.~{Wu},
  ``Optimal operation of soft open points in active distribution networks under
  three-phase unbalanced conditions,'' \emph{IEEE Trans. Smart Grid}, vol.~10,
  no.~1, pp. 380--391, 2019.

\bibitem{nizami2020multiagent}
M.~S.~H. {Nizami}, M.~J. {Hossain}, and E.~{Fernandez}, ``Multiagent-based
  transactive energy management systems for residential buildings with
  distributed energy resources,'' \emph{IEEE Trans. Ind. Informat.}, vol.~16,
  no.~3, pp. 1836--1847, 2020.

\bibitem{liu2020transactive}
Z.~{Liu}, L.~{Wang}, and L.~{Ma}, ``A transactive energy framework for
  coordinated energy management of networked microgrids with distributionally
  robust optimization,'' \emph{IEEE Trans. Power Syst.}, vol.~35, no.~1, pp.
  395--404, 2020.

\bibitem{yang2020transactive}
Z.~{Yang}, J.~{Hu}, X.~{Ai}, J.~{Wu}, and G.~{Yang}, ``Transactive energy
  supported economic operation for multi-energy complementary microgrids,''
  \emph{IEEE Transactions on Smart Grid}, 2020, to be published. DOI:
  \textcolor[rgb]{0,0,1}{10.1109/TSG.2020.3009670}.

\bibitem{yan2020distribution}
M.~{Yan}, M.~{Shahidehpour}, A.~{Paaso}, L.~{Zhang}, A.~{Alabdulwahab}, and
  A.~{Abusorrah}, ``Distribution network-constrained optimization of
  peer-to-peer transactive energy trading among multi-microgrids,'' \emph{IEEE
  Trans. Smart Grid}, 2020, to be published. DOI:
  \textcolor[rgb]{0,0,1}{10.1109/TSG.2020.3032889}.

\bibitem{cao2016benefits}
W.~Cao, J.~Wu, N.~Jenkins, C.~Wang, and T.~Green, ``Benefits analysis of soft
  open points for electrical distribution network operation,'' \emph{Appl.
  Energy}, vol. 165, pp. 36--47, 2016.

\bibitem{wang2015SNOP}
C.~{Wang}, C.~{Sun}, P.~{Li}, J.~{Wu}, F.~{Feng}, and Y.~{Yu}, ``{SNOP}-based
  operation optimization and analysis of distribution networks,'' \emph{Automa.
  Electric Power Syst.}, vol.~39, no.~9, pp. 82--87, 2015.

\bibitem{du2020Intelligent}
Y.~{Du} and F.~{Li}, ``Intelligent multi-microgrid energy management based on
  deep neural network and model-free reinforcement learning,'' \emph{IEEE
  Trans. Smart Grid}, vol.~11, no.~2, pp. 1066--1076, 2020.

\bibitem{li2017coordinated}
P.~{Li}, H.~{Ji}, C.~{Wang}, J.~{Zhao}, G.~{Song}, F.~{Ding}, and J.~{Wu},
  ``Coordinated control method of voltage and reactive power for active
  distribution networks based on soft open point,'' \emph{IEEE Trans. Sustain.
  Energy}, vol.~8, no.~4, pp. 1430--1442, 2017.

\bibitem{Oikonomou2019deliverable}
K.~{Oikonomou}, M.~{Parvania}, and R.~{Khatami}, ``Deliverable energy
  flexibility scheduling for active distribution networks,'' \emph{IEEE Trans.
  Smart Grid}, 2019, to be published. DOI:
  \textcolor[rgb]{0,0,1}{10.1109/TSG.2019.2927604}.

\bibitem{li2019optimal}
P.~Li, H.~Ji, C.~Wang, J.~Zhao, G.~Song, F.~Ding, and J.~Wu, ``Optimal
  operation of soft open points in active distribution networks under
  three-phase unbalanced conditions,'' \emph{IEEE Trans. Smart Grid}, vol.~10,
  no.~1, pp. 380--391, 2019.

\bibitem{liu2018energysharing}
N.~{Liu}, M.~{Cheng}, X.~{Yu}, J.~{Zhong}, and J.~{Lei}, ``Energy-sharing
  provider for {PV} prosumer clusters: A hybrid approach using stochastic
  programming and stackelberg game,'' \emph{IEEE Trans. Ind. Electron.},
  vol.~65, no.~8, pp. 6740--6750, 2018.

\bibitem{li2016sufficient}
Z.~{Li}, Q.~{Guo}, H.~{Sun}, and J.~{Wang}, ``Sufficient conditions for exact
  relaxation of complementarity constraints for storage-concerned economic
  dispatch,'' \emph{IEEE Trans. Power Syst.}, vol.~31, no.~2, pp. 1653--1654,
  2016.

\bibitem{yang2019real}
X.~Yang, Y.~Zhang, H.~He, S.~Ren, and G.~Weng, ``Real-time demand side
  management for a microgrid considering uncertainties,'' \emph{IEEE Trans.
  Smart Grid}, vol.~10, no.~3, pp. 3401--3414, 2019.

\bibitem{jiang2019stochastic}
Y.~Jiang, C.~Wan, J.~Wang, Y.~Song, and Z.~Y. Dong, ``Stochastic receding
  horizon control of active distribution networks with distributed
  renewables,'' \emph{IEEE Trans. Power Syst.}, vol.~34, no.~2, pp. 1325--1341,
  2019.

\bibitem{MultiFlexibility2020Marco}
M.~R.~M. {Cruz}, D.~Z. {Fitiwi}, S.~F. {Santos}, S.~J. P.~S. {Mariano}, and
  J.~P.~S. {Catalao}, ``Multi-flexibility option integration to cope with
  large-scale integration of renewables,'' \emph{IEEE Trans. Sustain. Energy},
  vol.~11, no.~1, pp. 48--60, 2020.

\bibitem{li2018energypricing}
Y.~Li, C.~Feng, F.~Wen, K.~Wang, and Y.~Huang, ``Energy pricing and management
  for park-level energy internets with electric vehicles and power-to-gas
  devices,'' \emph{Automa. Electric Power Syst.}, vol.~42, no.~16, pp. 1--10,
  2018.

\bibitem{yang2020flexibility}
X.~Yang, C.~Xu, H.~He, W.~Yao, J.~Wen, and Y.~Zhang, ``Flexibility provisions
  in active distribution networks with uncertainties,'' \emph{IEEE Trans.
  Sustain. Energy}, 2020, to be published. DOI:
  \textcolor[rgb]{0,0,1}{10.1109/TSTE.2020.3012416}.

\bibitem{yang2019event}
X.~Yang, Y.~Zhang, H.~Wu, and H.~He, ``An event-driven {ADR} approach for
  residential energy resources in microgrids with uncertainties,'' \emph{IEEE
  Trans. Ind. Electron.}, vol.~66, no.~7, pp. 5275--5288, 2019.

\end{thebibliography}
\bibliographystyle{IEEEtran}


\end{document}